\documentclass[review,onefignum,onetabnum]{siamart250211}
\PassOptionsToPackage{colorlinks=true,
            linkcolor=blue,    % for \ref, \cref, etc.
            citecolor=teal,    % for \cite
            urlcolor=magenta }{hyperref}

\usepackage[letterpaper, total={6in, 9in}]{geometry}
\usepackage{graphicx}
\usepackage{amsmath, amsfonts, latexsym}
\usepackage{enumitem}   
\usepackage{float}
\usepackage{url} % Insert URL addresses within the document
\usepackage[numbers]{natbib} % or plainnat
\usepackage{algorithm}
\usepackage{algorithmic}

\usepackage{cleveref}
\usepackage{tikz}
\usepackage{multirow}
\usepackage{pgfplots}
\usepackage{xifthen}
\usepackage{makecell}
\usepackage{subcaption}
\usepackage{caption}
\usepackage{amssymb}
\usepackage{booktabs}
\newcommand{\deleted}[1]{%
}
\usetikzlibrary{matrix}
\usepgfplotslibrary{groupplots}
\pgfplotsset{compat=newest}
\usetikzlibrary{shapes.geometric, arrows}

\tikzstyle{start} = [rectangle, rounded corners, 
minimum width=2cm, 
minimum height=1cm,
text centered, 
draw=black, 
text width=4cm, 
]

\tikzstyle{stop} = [rectangle, rounded corners, 
minimum width=2cm, 
minimum height=1cm,
text centered, 
draw=black, 
text width=1.5cm, 
]

\tikzstyle{io} = [trapezium, 
trapezium stretches=true, % A later addition
trapezium left angle=70, 
trapezium right angle=110, 
minimum width=3cm, 
minimum height=1cm, text centered, 
text width=4cm, 
draw=black, 
]

\tikzstyle{process} = [rectangle, 
minimum width=2cm, 
minimum height=0.7cm, 
text centered, 
text width=3cm, 
draw=black, 
]

\tikzstyle{decision} = [diamond, 
minimum width=2.5cm, 
minimum height=1cm, 
text centered, 
text width=2cm, 
draw=black, 
]
\tikzstyle{arrow} = [thick,->,>=stealth, line width=1.2pt]
\tikzset{line/.style={thick}}

\usetikzlibrary{arrows.meta}

\usepackage{lipsum}
\usepackage{amsfonts}
\usepackage{graphicx}
\usepackage{epstopdf}
\ifpdf
  \DeclareGraphicsExtensions{.eps,.pdf,.png,.jpg}
\else
  \DeclareGraphicsExtensions{.eps}
\fi

\newsiamremark{remark}{Remark}
\newsiamremark{hypothesis}{Hypothesis}
\crefname{hypothesis}{Hypothesis}{Hypotheses}
\newsiamthm{claim}{Claim}
\newsiamremark{fact}{Fact}
\crefname{fact}{Fact}{Facts}

\headers{An Adaptive Mixed Precision and Dynamically Scaled
PCG}{Y. Guo, E. de Sturler, and T. Warburton}

\title{An Adaptive Mixed-Precision and Dynamically Scaled Preconditioned Conjugate Gradient Algorithm
\thanks{Submitted to the editors 05/06/2025.
\funding{This work was supported under NSF DMS 2208470. The first author acknowledges the generous support from the John K. Costain Graduate Fellowship. The third author was supported in part by the John K. Costain Faculty Chair in Science at Virginia Tech. }}}

\author{Yichen Guo\thanks{Department of Mathematics,  Virginia Tech, 
Blacksburg, VA 24061
	 (\email{ycguo@vt.edu}, \email{sturler@vt.edu}, \email{tcew@vt.edu})} 
  \and Eric de Sturler\footnotemark[2] 
  \and Tim Warburton\footnotemark[2]
}

\usepackage{amsopn}

\definecolor{mylightyellow}{rgb}{1, 1, 0.9}
\definecolor{mylightblue}{rgb}{0.75, 0.85, 0.917}
\definecolor{myindigo}{rgb}{0.23, 0.18, 0.53}    % Slightly gray Indigo
\definecolor{mycyan}{rgb}{0.96, 0.74, 0.51}% Slightly gray Cyan
\definecolor{myyellow}{rgb}{0.96, 0.8, 0.3}% Slightly gray Cyan
\definecolor{myteal}{rgb}{ 0.3137,    0.5216,    0.5647}      % Slightly gray Teal
\definecolor{mygreen}{rgb}{0.13, 0.47, 0.27}     % Slightly gray Green
\definecolor{mylightgreen}{rgb}{0.74, 0.88, 0.64} % Light Green
\definecolor{myolive}{rgb}{0.57, 0.57, 0.33}     % Slightly gray Olive
\definecolor{mysand}{rgb}{0.8157,    0.8078,    0.6627}      % Slightly gray Sand
\definecolor{myrose}{rgb}{0.96, 0.65, 0.63}   % Slightly gray Rose
\definecolor{mywine}{rgb}{0.52, 0.25, 0.38}      % Slightly gray Wine
\definecolor{mypurple}{rgb}{0.3922,    0.1490,    0.4039}    % Slightly gray Purple
\definecolor{mylightpurple}{rgb}{0.74, 0.58, 0.78}     % Light Purple
\definecolor{mybrown}{rgb}{0.63, 0.39, 0.25}  
\definecolor{mydeepblue}{rgb}{0.2, 0.47, 0.75}    % Deep Blue
\definecolor{vtorange}{rgb}{0.8980,    0.4588,    0.1216}
\definecolor{vtgray}{rgb}{0.4588,    0.4706,    0.4824}
\definecolor{vtmaroon}{rgb}{0.5255,    0.1216,    0.2549}
\definecolor{lightgrey}{RGB}{214,214,214}
\definecolor{lightorange}{RGB}{255,169,66}

\def\eps{\varepsilon}

\def\C{\mathbb{C}}

\def\CNN{\C^{N\times N}\kern-2pt}

  \def\BA{{\bf A}} 
\def\Bb{{\bf b}}   
   
\def\Bd{{\bf d}}   
\def\Be{{\bf e}}

  \def\BI{{\bf I}}

  \def\BM{{\bf M}}

\def\Bp{{\bf p}}  \def\BP{{\bf P}} \def\CP{{\cal P}}
\def\Bq{{\bf q}}   
\def\Br{{\bf r}}   
   
  \def\BT{{\bf T}} 
   
  \def\BV{{\bf V}} 
   
\def\Bx{{\bf x}}   
\def\By{{\bf y}}   
\def\Bz{{\bf z}}

\def\Bzeta{\boldsymbol{\zeta}}
\def\Bxi{\boldsymbol{\xi}}
\def\Btheta{\boldsymbol{\theta}}
\def\Weta{\widehat{\eta}}

\def\epsS{\varepsilon_{\mathrm{32}}}
\def\epsD{\varepsilon_{\mathrm{64}}}
\def\BAs{\BA{\kern-1pt}}      % version of bold A for use in A^k
\def\BPs{\BP{\kern-1.2pt}}    % version of bold P for use in P_j
\def\BVs{\BV{\kern-1.2pt}}    % version of bold V for use in V_j
\def\CPs{\CP{\kern-1.2pt}}    % version of cal P for use in P_k

\def\tol{\textit{tol}}
\newcommand{\DPCG}{DS-PCG}
\newcommand{\APCG}{AMP-PCG}

\def\epsxk{\eps_{x,k}}

\def\epsqk{\eps_{q,k}}
\def\epsrk{\eps_{r,k}}
\def\uzk{u_{z,k}}
\def\urk{u_{r,k}}

\newtheorem{assumption}[theorem]{Assumption}

\newtheorem{example}{Example}

\crefname{hypothesis}{Hypothesis}{Hypotheses}
\def\lp{\left(}
\def\rp{\right)}

\def\Hz{\widehat{{\bf z}}}
\def\Hp{\widehat{{\bf p}}}
\def\Hq{\widehat{{\bf q}}}
\def\Hx{\widehat{{\bf x}}}

\def\Hr{\widehat{{\bf r}}}
\def\Hrho{\widehat{\rho}}
\def\Hbeta{\widehat{\beta}}
\def\Halpha{\widehat{\alpha}}
\def\Hgamma{\widehat{\gamma}}

\def\blue{\textcolor{blue}}
\def\orange{\textcolor{vtorange}}

\let\oldnl\nl
\newcommand{\nlnonumber}{\renewcommand{\nl}{\let\nl\oldnl}}

\newcommand{\exref}[1]{Ex.~\ref{#1}}

\usepackage{lineno}

\ifpdf
\hypersetup{
  pdftitle={An Adaptive Mixed Precision and Dynamically Scaled Preconditioned Conjugate Gradient Algorithm},
  pdfauthor={Yichen Guo, Eric de Sturler, Tim Warburton}
}
\fi

\begin{document}

\maketitle

% REQUIRED
\begin{abstract}
We propose an adaptive mixed precision and dynamically scaled preconditioned conjugate gradient algorithm (\APCG). It dynamically adjusts the precision  for storing vectors and computing, exploiting low precision  when appropriate, while maintaining a convergence rate and accuracy comparable to that of double precision PCG. 
Our mixed precision strategy consists of three main components:
(1)  The residual and matrix-vector product are initially computed in double precision, and the algorithm switches these to single precision based on the chosen convergence tolerance and an estimate of the residual gap. 
(2) Depending on the eigenvalue distribution, the preconditioned residual and search direction are either  in half precision throughout the iterations or initially in double precision and then stepwise reduced to single and half precision. 
(3) A dynamically scaled residual is used at every iteration to mitigate underflow in half precision. 
We provide theoretical support for our estimates and we demonstrate the effectiveness of \APCG{} through numerical experiments, highlighting both its robustness and the significant performance gains ($1.63\times$ speedup) achieved compared to  double precision PCG on a GPU.
\end{abstract}

% REQUIRED

% REQUIRED

\section{Introduction}

In this paper, we propose an {adaptive mixed precision and dynamically-scaled Preconditioned Conjugate Gradient (PCG)} algorithm, referred to as \APCG, for solving symmetric positive definite linear systems. During the PCG iterations, \APCG{} 
dynamically adapts the precision 
in which iteration vectors are stored and 
operations are carried out, using low-precision where appropriate, while ensuring convergence rate and accuracy comparable to  double
precision implementation.
This strategy can significantly reduce 
the runtime for the linear solver.
To mitigate the potential 
downsides of using 
mixed precision, our proposed algorithm 
satisfies three properties: 
\begin{enumerate}[label=\textbf{P\arabic*.}]
\item It ensures that the
true residual reaches the
chosen convergence tolerance by keeping the so-called
residual gap sufficiently small \cite{greenbaum1989behavior, greenbaum1997estimating}.
\label{item:p1}
\item It maintain roughly the rate 
of convergence that double 
precision PCG would achieve.
\label{item:p2}
\item It uses dynamic scaling 
of the iteration vectors, when using
FP16, to 
ensure that vectors (and vector components) can be represented 
within a relatively small
range.
\label{item:p3}
\end{enumerate}

In numerical simulations, solving large-scale linear systems is computationally intensive and often the most time-consuming part. 
Recent hardware accelerators provide significantly higher arithmetic throughput at lower precisions.
Moreover, using lower precision storage  reduces both memory usage and bandwidth demands, substantially accelerating data transfers, which in turn allows the algorithm to take better advantage of the low precision higher arithmetic throughput.
These trends motivate the development of algorithms that 
aggressively 
exploit 
reduced precision to improve the efficiency of sparse linear solvers.

\subsection{Literature review}
The use of multiple precisions in solving linear systems has been  explored. 
One well-studied approach is iterative refinement, where the solution is improved iteratively by solving a sequence of correction equations that account for the residual error. At each step, the residual is computed in high precision, while the correction equations are solved in lower precision using either a direct solver, often based on a low precision LU-factorization~\cite{langou2006exploiting, carson2018accelerating}, or Krylov subspace methods~\cite{amestoy2024five}. Mixed precision iterative refinement has been extensively studied in both theory and practice~\cite{goddeke2007performance, wilkinson2023rounding, higham2021exploiting, carson_mixed_2023}. A potential drawback of this approach is that the Krylov subspace is discarded upon each outer refinement iteration, potentially leading to a significantly higher total iteration count compared with a  double precision solve.

Beyond iterative refinement, various works have incorporated low precision arithmetic into specific components of Krylov subspace and other iterative methods~\cite{anzt2015adaptive}. 
For example, low precision preconditioners have demonstrated the ability to maintain convergence while accelerating computation~\cite{gobel2021mixed, kawai2022low, tian2024mixed}. 
In addition, matrix-vector products (matvec) can be performed in mixed precision, either statically or adaptively based on the magnitude of entries, to reduce data movement and computational cost~\cite{maddox2022low, graillat2024adaptive, yang2024mille}. 
Mixed precision restarted GMRES methods have also been proposed~\cite{lindquist2021accelerating}, where only the residual and solution updates are  in double precision, while other vectors are  in single precision.
Specifically for PCG, Clark et al.~\cite{clark2010solving} investigate the use of half precision multigrid preconditioning and single  precision arithmetic in lattice QCD solvers, where the updated residual produced by PCG is periodically replaced with the true residual computed in high precision.
In \cite{sleijpen1996reliable}, a strategy that determines when to perform this residual replacement without affecting the convergence is studied. 
Maddox et al.~\cite{maddox2022low} propose performing matvec operations in FP16 while using full orthogonalization to maintain the orthogonality of the residuals and applying log-sum-exp transformations to mitigate rounding and overflow errors. 

Moreover, our work is closely related to the theory of finite precision CG. 
In such settings, rounding errors due to the finite precision arithmetic leads to two primary effects: a reduction in the attainable accuracy and a delay in convergence, as discussed in \cite{paige1980accuracy,van2000residual,greenbaum1997estimating, gutknecht2000accuracy}. 
These challenges are magnified in communication-avoiding variants, such as pipelined CG and $s$-step CG, 
which have larger local rounding errors relative to CG;  see \cite{carson2018numerical, cools2018analyzing, carson2018adaptive, carson2022mixed} for detailed analysis of the resulting limitations in accuracy and convergence.

Mixed precision Krylov methods can also be analyzed as \emph{inexact Krylov subspace methods}. In~\cite{golub1999inexact}, the perturbations introduced by preconditioning using inner solves is analyzed. In~\cite{bouras2000relaxation, bouras2005inexact, du2008varying}, the matvec is computed approximately using an inner iteration, and adaptive stopping criteria are proposed for inner iteration to ensure convergence of the outer iteration. Their analysis shows that the allowable inaccuracy in the matvec can increase as the outer iteration progresses toward convergence.

\subsection{\APCG{}}

We distinguish two choices of precision in \APCG{}: 
one for the matvec and residual vector, which mainly controls the final  attainable residual accuracy \cite{greenbaum1997estimating},
and the second precision for the preconditioned residual and search direction vectors, which primarily affects the convergence rate \cite{golub1999inexact, notay2000flexible, greenbaum2021convergence}.
It has been shown in \cite{greenbaum1997estimating} that roundoff errors in residual updates can limit the maximum attainable residual accuracy. Based on this analysis, we propose an attainable accuracy indicator and a precision switching criterion  to reduce the precision of residual vector and matvec 
from FP64 to FP32. 
The second choice of precision concerns the preconditioned residual and the search direction vector.
If the system matrix is well-conditioned or equipped with an effective preconditioner, 
the convergence is typically linear, and 
the preconditioned residual and search direction can be stored in FP16 from the start, while still achieving convergence rates comparable to PCG in FP64.
This is the case with, for example, $p$-multigrid when applied to high-order finite element methods (FEM) for Poisson equations \cite{karakus2019gpu, min2022optimization}.
We demonstrate the convergence of our algorithm with a multigrid preconditioner and search direction vectors in half precision in \cref{sec: well conditioned matrices}. 
For linear systems without many large outlying eigenvalues and 
if the convergence is not linear,
we begin using FP64 and lower the precision of the preconditioned residual and search direction to FP32 and FP16 according to the precision selector in \APCG{}. In these cases,  convergence remains close to that of PCG in FP64. For matrices with large, outlying eigenvalues, using vectors in low precision may lead to frequent recurrence  of the eigenvectors corresponding to the large(st) eigenvalues in the construction of Krylov space~\cite{greenbaum1992predicting}, which can degrade convergence. In such cases, the algorithm adapts the precision more conservatively.

To ensure that all vectors remain within the representable range of the current precision, 
\APCG{}  normalizes the residual before preconditioning (other scaling factors can be used). 
It can be shown that the dynamically scaled PCG is equivalent to PCG in exact arithmetic. 
This strategy is similar to scaling in iterative refinement~\cite{higham2021exploiting}, but is applied in each iteration rather than applied in the outer loop. Scaling in low precision has also been used effectively in deep learning~\cite{pytorch_amp} and fluid dynamics simulations~\cite{klower2022fluid} to prevent underflow and overflow. 
In summary, the \APCG{} algorithm follows almost the same steps as PCG but employs mixed precision for different vectors adaptively, and it typically achieves convergence and accuracy comparable to a double precision PCG while reducing computational cost and data movement.

In the following sections, we first motivate the careful design of mixed precision PCG with an illustrative example in \cref{sec:motivation}. Then, in \cref{sec: theory}, we review theoretical results of finite precision CG for attainable accuracy and inexact Krylov methods for convergence rates. 
In \cref{sec:dspcg}, we introduce a dynamically scaled PCG method designed to prevent overflow and underflow in low precision arithmetic. 
We present our main algorithm, adaptive mixed precision and dynamically scaled PCG, and precision switch criteria in \cref{sec: mixed precision}. 
Finally, in \cref{sec:numerical}, we evaluate the method’s performance by comparing its convergence behavior to that of double precision PCG on various linear systems and assessing runtime speedups in a large-scale GPU-accelerated high-order finite element solver.

\section{Motivating example}
\label{sec:motivation}

We consider the solution of  the $n \times n$ linear system
\begin{equation}
    \BA \Bx = \Bb,
    \label{eq:Axb}
\end{equation}
where $\BA$ is a real, symmetric, positive definite (SPD) matrix.
The preconditioned conjugate gradient (PCG) method \cite{hestenes1952methods} is widely used for solving large-scale, sparse, SPD systems. An outline of the PCG algorithm is given in \cref{alg:PCG}.

In exact arithmetic, the updated residual $\Br_{k+1}$ exactly equals the true residual $\Bb - \BA\Bx_{k+1}$, the search directions $\Bp_k$ are $\BA$-orthogonal to all previous search directions, and the algorithm converges to the exact solution in at most $n$ iterations. 
However, in finite precision arithmetic, rounding errors disrupt these properties \cite{greenbaum1992predicting, meurant2006lanczos}. 
As a result,
the updated residual $\Br_{k+1}$, obtained from line~13 of \cref{alg:PCG},  may deviate from the true residual $\Bb - \BA\Bx_{k+1}$, due to the rounding error in line~12 and 13.
Additionally, the search direction $\Bp_k$ may lose its global $\BA$-orthogonality with previous directions, maintaining only local orthogonality, which can result in a degraded convergence rate.  
Consequently, the algorithm may require more than \(n\) iterations to converge.
We present a motivating example to demonstrate how finite precision affects the convergence of PCG.
\begin{example}
    Let  $\BA \in \mathbb{R}^{100 \times 100}$ have eigenvalues $\lambda_1, \dots, \lambda_{100}$,  defined as follows:
\begin{equation}
\begin{aligned}
& \lambda_{100} = 1, \quad 
\lambda_i = 10^{-3} + \frac{i-1}{99} \left( \frac{17}{20} \right)^{100 - i}, \quad \text{for } i = 95, \dots, 99,\\
& [\lambda_1, \dots, \lambda_{95}] = \texttt{linspace}(10^{-3}, \lambda_{95},\, 95).  
\end{aligned}
\label{eq:motivating A}
\end{equation}
Here, \texttt{linspace}$(a, b, N)$ denotes $N$ evenly spaced values between $a$ and $b$.
The matrix $\BA$ contains 5 isolated large eigenvalues and 95 evenly spaced  eigenvalues. 
The eigenvector matrix of $\BA$ is  a random orthonormal matrix. 
The right-hand side  $\Bb$ is a vector of all ones. The preconditioner is the identity matrix, and the initial guess $\Bx_0$ is the zero vector. The stopping tolerance is $\tol = 10^{-10}$.
\label{ex: motivating}
\end{example}

\begin{figure}[h!]
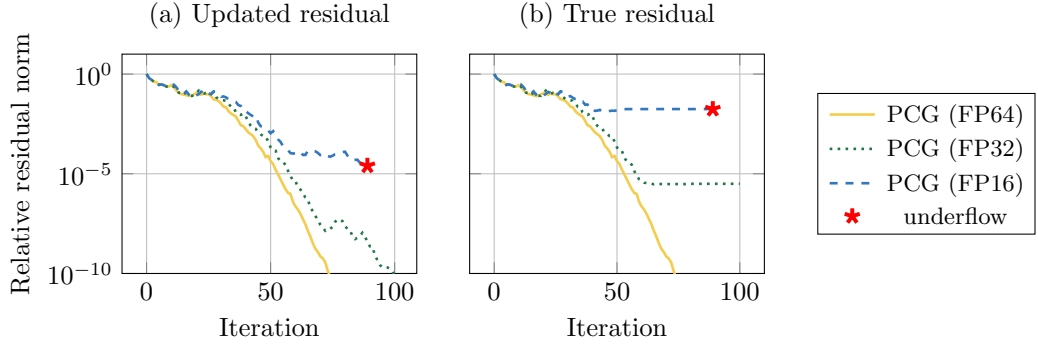

\centering
\begin{tikzpicture}
\begin{groupplot}[
    group style={
        group name=myplot,
        group size=2 by 1,
        horizontal sep=20pt,
        vertical sep=2cm
    },
    height=4.5cm,
    width=5.5cm
]
  % First plot with y-axis labels
  \nextgroupplot[
  title={(a) Updated residual},
    xlabel={Iteration},
    ylabel={Relative residual norm},
    ymode=log,
    ymin=1e-10,
    ymax=10,
     grid=both,
         ytick={1e0, 1e-5, 1e-10, 1e-14},
  ]
  \input{2.Motivation/data/motivating_init_1}
  
  % Second plot without y-axis tick labels
  \nextgroupplot[
    title={(b) True residual},
    xlabel={Iteration},
    ymode=log,
    ymin=1e-10,
    ymax=10,
    yticklabels={},
        ytick={1e0, 1e-5, 1e-10, 1e-14},
     grid=both,
  ]
  \input{2.Motivation/data/motivating_init_2}
 
\end{groupplot}

\path (myplot c1r1.south west|-current bounding box.south)--
      coordinate(legendpos)
      (myplot c2r1.south east|-current bounding box.south);

\matrix[
    matrix of nodes,
    anchor=west, % align the left edge of the legend with the anchor point
    draw,
    inner sep=0.2em,
    nodes={font=\small}
] at ([xshift=2em]myplot c2r1.east) % shift a bit to the right
{
    \ref{plot:fp64 motivation} & PCG (FP64) \\
    \ref{plot:fp32} & PCG (FP32) \\
    \ref{plot:fp16 motivation} & PCG (FP16) \\
    \ref{plot:fp16underflow} & underflow 
    \\
};
\end{tikzpicture}
\caption{
Convergence of (a) the updated residual, and (b) the true residual, computed as $\Bb - \BA\Bx_k$, for PCG implemented  in half, single, and double precision for \cref{ex: motivating}. 
The matrix $\BA$ is defined in \cref{eq:motivating A}. 
Compared with FP64, both FP16 and FP32 yield lower attainable accuracy and slower convergence.
In the FP16 case, $\rho_k$ eventually underflows to zero due to limited dynamic range, leading to breakdown, indicated by the red star. 
This example highlights that naive use of low precision in PCG  can lead to  reduced accuracy, delayed convergence, and data underflow.
}
\label{fig: motivation example}
\end{figure}

We solve the linear system using the PCG algorithm in  half, in single, and in double precision. Figure~\ref{fig: motivation example} shows the convergence  of  the relative updated residual norm and the relative true residual norm for these three implementations. Here, the updated residual refers to the residual $\Br_k$ computed in line~13 of Algorithm~\ref{alg:PCG}, while the true residual is given by $\Bb - \BA\Bx_k$.
Initially, the convergence of PCG in half or single precision closely matches that of the double precision implementation. However, after approximately 30 iterations, the half and double precision PCG  diverge.
As the residual norm continues to decrease, in half precision, $\Bz_k$ eventually underflows to zero due to the limited dynamic range, causing $\rho_k$ to be zero and consequently  $\beta_k$ evaluation induces NaNs. This underflow is marked with a red star in  \cref{fig: motivation example}.
Moreover, due to the reduced numerical accuracy in low precision formats, the updated residual norm deviates from the true residual norm in both half and single precision implementations. 
As a result, the final accuracy in low precision is substantially lower than that achieved in double precision.
Additionally, the convergence of PCG in both half and single precision is slightly slower than in double precision after the 25-th iteration. 
This example illustrates three key challenges  using PCG in low precision:  unacceptable accuracy, degraded convergence rate, and data underflow.

Our goal is to develop an algorithm that satisfies \ref{item:p1}, \ref{item:p2}, \ref{item:p3},
while using low precision in selected steps and variables of  PCG, thereby reducing computational cost and data movement. 
In the following sections, we review key results from finite precision sensitivity analysis  of PCG and inexact Krylov methods to understand what affects the attainable accuracy and the rate of convergence of PCG.
Additionally, we introduce a dynamically scaling to mitigate data underflow.
We propose an algorithm that selects appropriate precision levels for $\Bz_k$, $\Bp_k$, $\Br_k$, and $\Bq_k$ in \cref{alg:MPDS-PCG}.
The motivating example is revisited in \cref{sec: revisit motivating}.

\section{Analysis}
\label{sec: theory}

We recall the finite precision analysis on the attainable accuracy and the theory of inexact Krylov subspace methods about the convergence rate of PCG. To address potential data underflow, we introduce the dynamically scaled PCG algorithm (\DPCG). The analysis and the \DPCG{} algorithm provide the foundation for the adaptive mixed precision and dynamically scaled PCG method described in \cref{sec: mixed precision}.

\subsection{(P1) Attainable accuracy}
\label{sec:attainable accuracy}

For the attainable accuracy of PCG, we largely follow 
\cite{greenbaum1997estimating}, with modifications tailored to our specific context.
This analysis  holds for any coupled recurrence of the type
\begin{eqnarray}
\label{eq:basic_rec_1}
\Bx_k & = & \Bx_{k-1}+\alpha_{k-1}\Bp_{k-1}, \\
\label{eq:basic_rec_2}
\Br_k & = & \Br_{k-1} - \alpha_{k-1}\BA\Bp_{k-1} .
\end{eqnarray}
Hence, the analysis remains valid even if the mixed precision PCG deviates substantially from PCG in exact arithmetic and also under changes in other parts of the algorithm.

Let $\epsxk$, $\epsrk$, and $\epsqk$ denote the roundoff errors in computing $\Bx_k$, $\Br_k$, and $\Bq_k$, respectively, as specified in lines 12, 13, and 9 of \cref{alg:PCG}.
A finite precision implementation leads to the following (extended) recurrence (variables denote the finite precision computed variables),
\begin{align}
\Bx_k & =  \Bx_{k-1}+\alpha_{k-1}\Bp_{k-1} + \Bxi_{k}, \\
\Br_k & =  \Br_{k-1} - \alpha_{k-1} \BA\Bp_{k-1} + \Btheta_k,
\label{eq:res perturb}
\end{align}
where
\begin{align}
  \|\Bxi_k\| & \leq 
  \epsxk \|\Bx_{k-1}\| +(2\epsxk + \epsxk^2)\|\alpha_{k-1}\Bp_{k-1}\|, \label{eq:perturb x} \\
  \|\Btheta_k\| & \leq 
  \epsrk\| \Br_{k-1}\| +
  (2\epsrk + \epsrk^2)\|\alpha_{k-1} \BA \Bp_{k-1}\| +
  (1+2\epsrk + \epsrk\epsqk)\|\alpha_{k-1}\Bd_{k-1}\|.
\label{eq:pertub bound}
\end{align}
These bounds follow directly from \cite{greenbaum1997estimating} using the  specific roundoff errors above.
The vector $\Bd_{k-1}$ represents the numerical error arising from the floating point evaluation of $\BA\Bp_{k-1}$, given by,
$\Bd_{k-1} =  \mathrm{fl}(\BA\Bp_{k-1}, \epsqk) - \BA\Bp_{k-1}$, where $ \mathrm{fl}(\cdot, \epsqk)$ denotes the result of floating point arithmetic  with $\epsqk$. 
We define $\Bzeta_k$  as the difference between the true residual $\Bb - \BA\Bx_k$ and the updated residual $\Br_k$,
\begin{equation}
    \Bzeta_k = \Bb - \BA\Bx_k - \Br_k = \lp \Bb - \BA\Bx_0 - \Br_0\rp -  \sum\limits_{t=1}^{k} \lp \BA \Bxi_t + \Btheta_t \rp.
    \label{eq:res_gap}
\end{equation}
Therefore, the true residual norm  $\|\Bb - \BA\Bx_k\|$ can be bounded as follows: 
\begin{equation}
   \left|  \|\Bzeta_k\| - \|\Br_k\| \right|  \leq  \|\Bb - \BA\Bx_k\| \leq   \|\Bzeta_k\| + \|\Br_k\|.
\end{equation}
In general, $\Br_k$, the updated residual, satisfies \cite{greenbaum1989behavior, greenbaum1997estimating},
\begin{eqnarray}
  \label{eq:comp_res_limit}
  \Br_k \rightarrow 0, &&  
  \quad \textrm{ for } k \rightarrow \infty.
\end{eqnarray}
When $\|\Bzeta_k\| \gg \|\Br_k\|$, the best achievable accuracy for $\|\Bb - \BA\Bx_k\|$ is approximately $\|\Bzeta_k\|$.
Therefore, to ensure that the true residual norm $\|\Bb - \BA\Bx_k\|$  converges to the chosen tolerance $\tol$ \footnote{for simplicity, we choose the stopping criterion based on the relative residual norm.},
the relative convergence tolerance should be set to 
\begin{eqnarray}
  \label{eq:comp_res_conv}
  { \|\Bzeta_k\| + \|\Br_k\|}
  \leq \tol\, {\|\Bb\|}.
\end{eqnarray}
Then, when this bound is satisfied at, say, iteration $m$,
we have 
\begin{eqnarray}
\label{eq:conv_true_res}
  \|\Bb - \BA\Bx_m\| \leq \|\Bzeta_m\|   + \|\Br_m\| 
  & \leq & \tol\,\|\Bb\|.
\end{eqnarray}
This approach assumes that
the {\em residual gap} \eqref{eq:res_gap} (or its upper bound),
stays below the convergence tolerance. 
If the {residual gap}  becomes too large, we can use several 
{\em residual replacement strategies} to enforce a sufficiently small gap while maintaining good convergence for the PCG algorithm \cite{van2000residual, sleijpen1996reliable, clark2010solving, carson2018numerical}. However, here we focus on controlling the gap by switching precision for various vectors at appropriate points in the algorithm. 
This strategy ensures that the attainable precision can be reached; however, it does not address the rate of convergence, which we consider next.

\subsection{(P2) Rate of convergence}
\label{sec: z and p precision}

Next, we consider how to maintain (close to) the rate of convergence achieved by PCG in double precision while exploiting mixed precision implementations. 
We follow the approach of \cite{golub1999inexact}, which considers the tolerance for an iteratively applied preconditioner. 
We can apply their analysis, for the purpose of this paper, 
to guide the precision of the preconditioner $\BM$  and the preconditioned residual  $\Bz_k$.
Moreover, as the quality of the new search direction, $\Bp_k = \Bz_k + \beta_{k-1}\Bp_{k-1}$, depends largely on $\Bz_k$ (the direction expanding the current Krylov subspace), there is no need to compute $\Bp_k$ in higher precision than $\Bz_k$. So, we use the same precision for both.   

The precision of $\Bz_k$  can be considered to reflect the accuracy with which the preconditioner has been computed. Following \cite{golub1999inexact}, we replace $\BM\Bz_k=\Br_k$  by 
$$\BM \Bz_k = \Br_k + \Be_k,$$
where $\Be_k$ is the perturbation in computing $\Bz_k$.  Theorem 3.6 in \cite{golub1999inexact} provides a bound for the rate of linear convergence for the Inexact Preconditioned CG (IPCG):
\begin{theorem}[Theorem 3.6 in \cite{golub1999inexact}]
Let 
$
 \delta \geq \left\|\Be_{k}\right\|_{\BM^{-1}}/\left\|\Br_{k}\right\|_{\BM^{-1}}.
$
If $ \delta'= 2 \sin \theta \sqrt{\kappa} < 1$ with $\theta = \arcsin  \delta$, $\kappa = \mathrm{cond}(\BM^{-1/2} \BA \BM^{-1/2})$, then IPCG converges, and for even $k$,
\begin{equation}
    \frac{\| \Br_{k+1} \|_{\BA^{-1}}}{\| \Br_1 \|_{\BA^{-1}}} \leq (\sigma K)^k,
    \label{eq:ipcg-convergence}
\end{equation}
where
\[
\sigma = \frac{\sqrt{\kappa'} - 1}{\sqrt{\kappa'} + 1}, \quad
K = \sqrt{ \frac{2 + 16  \delta' \kappa' / (\kappa' - 1)^2}{1 + \sigma^4} }, \quad
\text{and} \quad
\kappa' = \kappa \frac{1 +  \delta}{1 -  \delta}.
\]
\label{thm: 3.1}
\end{theorem}
In \cite{golub1999inexact}, the perturbation $ \delta$  is of the order $10^{-2}$. However, in the mixed precision context, using single precision or half precision for $\Bz_k$ and $\Bp_k$ usually leads to $ \delta \approx 10^{-7}$ and  $ \delta \approx 10^{-3}$, respectively. 
When $ \delta\ll 1$, the modified condition number $\kappa'$ remains close to $\kappa$. For example, if $\BM = \BI$ and $\Bz_k$ is computed in half precision, $ \delta \approx 10^{-3}$. Assuming $\kappa = 10^3$, we obtain $K \approx 1.06$, and the relative increase in $\kappa'$ compared to $\kappa$ is only about $0.2\%$. 
Thus, we anticipate if PCG is in a linear convergence regime, then using half precision for $\Bz_k$ and $\Bp_k$ does not significantly degrade the convergence rate.

\subsection{(P3) Dynamic Scaling}
\label{sec:dspcg}

As discussed above, we can use low precision vectors in PCG without significantly affecting convergence rate or attainable accuracy under certain conditions. 
Beyond convergence rate and accuracy, another critical issue in low precision is the risk of data underflow or overflow. The motivating example in \cref{sec:motivation} shows that underflow may occur 
when entries of the residual become small, causing the preconditioned residual to drop below the dynamic range of half precision.
To address this,  all vectors stored in half precision must be representable within a limited dynamic range. This can generally be achieved by consistently scaling all vectors during each PCG iteration. 
We introduce the dynamically scaled PCG algorithm (\DPCG), presented in \cref{alg:dspcg}, where the modifications relative to  PCG are highlighted in blue.
\noindent
\begin{minipage}[h]{0.48\textwidth}
\begin{algorithm}[H]
\caption{Preconditioned Conjugate Gradient (PCG) Algorithm}
\label{alg:PCG}
\begin{algorithmic}[1]
\REQUIRE SPD matrix $\mathbf{A}$, preconditioner $\mathbf{M}$, right-hand side $\mathbf{b}$, initial guess $\mathbf{x}_0$, stopping tolerance \tol, maximum number of iterations $m$
\ENSURE Approximate solution $\mathbf{x}$
\STATE Initialize: $\mathbf{r}_0 \gets \mathbf{b} - \mathbf{A} \mathbf{x}_0$, $\beta_{-1} \gets 0$
\FOR{$k = 0$ to $m-1$}
    \STATE $ \mathbf{z}_k \gets \mathbf{M}^{-1}\mathbf{r}_k$
    \STATE $\rho_k \gets \mathbf{r}_k^T \mathbf{z}_k$
    \IF{$k > 0$}
        \STATE $\beta_{k-1} \gets \rho_k / \rho_{k-1}$
    \ENDIF
    \STATE $\mathbf{p}_k \gets \mathbf{z}_k + \beta_{k-1} \mathbf{p}_{k-1}$ 
    \STATE $\mathbf{q}_k \gets \mathbf{A} \mathbf{p}_k$
    \STATE $\gamma_k \gets \mathbf{q}_k^T \mathbf{p}_k$
    \STATE $\alpha_k \gets \rho_k / \gamma_k$
    \STATE $\mathbf{x}_{k+1} \gets \mathbf{x}_k + \alpha_k \mathbf{p}_k$ \hfill 
    \STATE $\mathbf{r}_{k+1} \gets \mathbf{r}_k - \alpha_k \mathbf{q}_k$ \hfill
    \IF{$\|\mathbf{r}_{k+1}\| \leq \tol \,\|\Bb\|$}
        \STATE \textbf{break} 
    \ENDIF
\ENDFOR
\STATE $\mathbf{x} \gets \mathbf{x}_{k+1}$
\end{algorithmic}
\end{algorithm}
\end{minipage}%
\hfill
\begin{minipage}[h]{0.48\textwidth}
\begin{algorithm}[H]
\caption{Dynamically Scaled PCG (\DPCG)}
\label{alg:dspcg}
\begin{algorithmic}[1]
\REQUIRE SPD matrix $\mathbf{A}$, preconditioner $\mathbf{M}$, right-hand side $\mathbf{b}$, initial guess $\mathbf{x}_0$, stopping tolerance \tol, maximum number of iterations $m$, and scaling factors $\{\omega_k\}_{k=0}^{m}$
\ENSURE Approximate solution $\mathbf{x}$
\STATE Initialize: $\mathbf{r}_0 \gets \mathbf{b} - \mathbf{A} \mathbf{x}_0$, $\beta_{-1} \gets 0$
\FOR{$k = 0$ to $m-1$}
    \STATE \textcolor{blue}{$\mathbf{y}_k \gets \omega_k\,\mathbf{r}_k, \, \mathbf{z}_k \gets \mathbf{M}^{-1}\mathbf{y}_k$} 
    \STATE $\rho_k \gets \mathbf{r}_k^T \mathbf{z}_k$
    \IF{$k > 0$}
        \STATE $\beta_{k-1} \gets \rho_k / \rho_{k-1}$
    \ENDIF
    \STATE $\mathbf{p}_k \gets \mathbf{z}_k + \beta_{k-1} \mathbf{p}_{k-1}$
    \STATE $\mathbf{q}_k \gets \mathbf{A} \mathbf{p}_k$
    \STATE $\gamma_k \gets \mathbf{q}_k^T \mathbf{p}_k$
    \STATE $\alpha_k \gets \rho_k / \gamma_k$
    \STATE $\mathbf{x}_{k+1} \gets \mathbf{x}_k + \alpha_k \mathbf{p}_k$
    \STATE $\mathbf{r}_{k+1} \gets \mathbf{r}_k - \alpha_k \mathbf{q}_k$
    \IF{$\|\mathbf{r}_{k+1}\| \leq \tol\, \|\Bb\|$}
        \STATE \textbf{break}
    \ENDIF
\ENDFOR
\STATE $\mathbf{x} \gets \mathbf{x}_{k+1}$
\end{algorithmic}
\end{algorithm}
\end{minipage}
\vspace{0.5cm}

The key modification is to scale the residual vector prior to preconditioning (line 3 of \cref{alg:dspcg}) to ensure that $\|\By_k\|$ remains within a moderate range, 
thereby avoiding excessively small $\Bz_k$  values. The scaling factor $\omega_k$ in practice may be chosen as $\|\Br_k\|^{-1}$ or $\|\Bz_{k-1}\|^{-1}$. 
The scaling can be fused with preconditioning step to minimize overhead. 
Similar scaling strategies have been successfully used in other contexts involving mixed precision computation \cite{higham2019squeezing, klower2022fluid, pytorch_amp}.
For any sequence \(\{\omega_k\}_{k=0}^{m}\) with \(\omega_k \in \mathbb{R} \setminus \{0\}\) in \DPCG, the approximate solutions and residuals produced by PCG are exactly the same as those produced by \DPCG{} in exact arithmetic. This equivalence is formalized in the following theorem.

\begin{theorem}
Let $\Bx_k$ and $\Br_k$ denote the approximate solutions and residuals generated by PCG (\cref{alg:PCG}), and let $\Hx_k$ and $\Hr_k$  be those generated by  \DPCG{}(\cref{alg:dspcg}).
Given the same initial guesses $\Bx_0$ and a sequence of nonzero real scaling factors $\{\omega_k\}_{k=0}^{m}$, \DPCG{} and PCG iterates are identical in exact arithmetic for all $k$, that is,
\[
\Hx_k = \Bx_k \quad \text{and} \quad \Hr_k = \Br_k.
\]
\label{thm: dspcg}
\end{theorem}
The proof of \cref{thm: dspcg} is provided in \cref{apx: dspcg}.

\section{Adaptive mixed precision and dynamically scaled PCG algorithm}
\label{sec: mixed precision}

Based on the theoretical discussion in \cref{sec: theory}, we propose an adaptive mixed precision PCG algorithm with dynamic scaling that satisfies properties \ref{item:p1}, \ref{item:p2}, and \ref{item:p3}, while 
lowering the precision as much as possible to reduce computational cost and data movement.

Since the residual gap \cref{eq:res_gap} is influenced primarily by rounding errors in the computation of $\Bx_k$, $\Bq_k$, and $\Br_k$, while the rate of convergence depends primarily on rounding errors in computing  $\Bp_k$ and $\Bz_k$, 
we define two separate precision controls:  
\begin{enumerate}
    \item \( u_{r,k} \): the precision used to compute  \( \By_k \), and to compute and store \( \Bq_k \) and \( \Br_k \),
    \item \( u_{z,k} \): the precision used to  store \( \By_k \), and to  compute and store \( \Bz_k \) and \( \Bp_k \).
\end{enumerate}
The approximate solution \( \Bx_k \), all inner products, norms, and scalar-only operations are computed in double precision to ensure the accuracy of the computed solution.

To guarantee that we reach the target accuracy, we  
control $\|\Bzeta_k\|$ \cref{eq:res_gap}  by adaptively reducing \( u_{r,k} \) from double to single precision, based on a precision switch indicator  $\eta_k$ detailed in \cref{sec: criterion rq}.
The analysis of IPCG in \cref{sec: z and p precision} suggests that, within a linear convergence regime and under small perturbations in 
$\Bz_k$ and $\Bp_k$, the convergence rate remains comparable to that of double precision PCG.
Often,
when highly effective preconditioners (like multigrid) are used, the iteration exhibits  linear convergence, allowing  \( u_{z,k} \) to be FP16 throughout.
When linear convergence is not evident, we introduce two heuristic thresholds, 
\( \tau_{z,s} \) and \( \tau_{z,h} \), to estimate the onset of the linear regime.
These thresholds determine when $u_{z,k}$ should be reduced from double to single precision and from single to half precision, respectively. For our examples, we determine the thresholds experimentally. 

Let \( u_0 \) denote the initial (maximum) precision for $\uzk$, and let \( \tol \) be the stopping tolerance. We define the relative residual norm at iteration $k$ as $\nu_k = \frac{\|\Br_k\|}{\|\Bb\|}$.
The precision selector function below determines  \( u_{z,k} \) and \( u_{r,k} \) at iteration $k$: 
\begin{equation}
\begin{aligned}
    (u_{z,k}, u_{r,k}) & = \text{precisionSelector}\Big( \nu_k,\; u_0,\; \tau_{z,s},\; \tau_{z,h},\; \eta_k \Big), \quad \text{with}
    \\
u_{z,k} & =
\begin{cases}
u_0, & \text{if } \nu_k \ge \tau_{z,s}, \\[1mm]
\min\left\{u_0, \,\mathrm{FP32} \right\}, & \text{if } \tau_{z,h} \le \nu_k < \tau_{z,s}, \\[1mm]
\min\left\{u_0, \, \mathrm{FP16} \right\}, & \text{if } \nu_k < \tau_{z,h},
\end{cases}
\\
u_{r,k} & =
\begin{cases}
\mathrm{FP64}, \quad & \text{if } \eta_k \ge \tol, \\[1mm]
\mathrm{FP32}, \quad & \text{if } \eta_k < \tol,
\end{cases}
\end{aligned}
\label{eq:precision_selector}
\end{equation}
Here, \( \min(p_1, p_2) \) returns the lower precision format, and \( \eta_k \) is defined by \cref{eq: eta_k} or, in the case of linear convergence, by \cref{eq:linear criterion}.

The \emph{Adaptive Mixed Precision and Dynamically Scaled PCG} algorithm (\APCG), presented in \cref{alg:MPDS-PCG}, implements
the proposed algorithmic changes. 
Vectors highlighted in blue and orange reflect the precisions $\uzk$ and $\urk$, respectively.
\begin{algorithm}[htpb]
\caption{Adaptive Mixed Precision and Dynamically Scaled PCG with Initial Precision $u_0$ (\APCG($u_0$))}
\label{alg:MPDS-PCG}
\begin{algorithmic}[1]
\REQUIRE SPD matrix $\mathbf{A}$, preconditioner $\mathbf{M}$, right-hand side $\mathbf{b}$, initial guess $\mathbf{x}_0$,  
stopping tolerance $\tol$, precision switch thresholds $\tau_{z,h}$ and $\tau_{z,s}$, precision switch indicator $\eta_k$, 
maximum iteration count $m$, the default precision $u_0$, and a precision selector function
\ENSURE Approximate solution $\mathbf{x}$

\STATE $\mathbf{r}_0 \gets \mathbf{b} - \mathbf{A} \mathbf{x}_0$, $\beta_{-1} \gets 0$, $\delta_0 \gets \|\mathbf{r}_0\|$

\FOR{$k = 0$ to $m-1$}
    \STATE {Calculate $\eta_k$} \hfill //{Precision switch indicator for $u_{r,k}$}

    \STATE {\blue{$(u_{z,k}$}, \orange{$u_{r,k}$}) $\gets \texttt{precisionSelector} \left( \frac{\delta_k}{\|\mathbf{b}\|}, u_0, \tau_{z,s}, \tau_{z,h}, \eta_k \right)$} \hfill //{Select precisions following \cref{eq:precision_selector}}
    
    \STATE \blue{${\mathbf{y}}_k$}$ \gets $\orange{$\mathbf{r}_k$}$ / \delta_k$ 
    % \hfill //{Computed in $u_{r,k}$, stored in $u_{z,k}$}
    
    \STATE \textcolor{blue}{ $ {\mathbf{z}}_k \gets \mathbf{M}^{-1} {\mathbf{y}}_k$} 
    % \hfill //{Preconditioned residual in $u_{z,k}$}
    
    \STATE $\rho_k \gets $ \orange{$\mathbf{r}_k^T $}\blue{${\mathbf{z}}_k$} 
    % \hfill //{Computed in FP64}
    
    \IF{$k > 0$}
        \STATE $\beta_{k-1} \gets \rho_k / \rho_{k-1}$ 
        % \hfill //{Computed in FP64}
    \ENDIF
    
    \STATE \blue{${\mathbf{p}}_k \gets {\mathbf{z}}_k + \beta_{k-1} {\mathbf{p}}_{k-1}$} 
    % \hfill //{Search direction update in $u_{z,k}$}

    \STATE \orange{$\mathbf{q}_k \gets \mathbf{A}$}\blue{$ {\mathbf{p}}_k$ }
    % \hfill //{Matrix-vector product in $u_{r,k}$}
    
    \STATE $\gamma_k \gets$\orange{$\mathbf{q}_k^T$}\blue{$ {\mathbf{p}}_k$} 
    % \hfill //{Computed in FP64}
    
    \STATE $\alpha_k \gets \rho_k / \gamma_k$ 
    % \hfill //{Computed in FP64}
    
    \STATE $\mathbf{x}_{k+1} \gets \mathbf{x}_k + \alpha_k $\blue{${\mathbf{p}}_k$} \hfill //{Update solution in FP64}
    
    \STATE \orange{$\mathbf{r}_{k+1} \gets \mathbf{r}_k - \alpha_k \mathbf{q}_k$} \hfill 
    % //{Update residual in $u_{r,k}$}
    
    \STATE $\delta_{k+1} \gets $\orange{$\|\mathbf{r}_{k+1}\|$}

    \IF{$\delta_{k+1} \leq \tol \, \|\Bb\|$}
        \STATE \textbf{break}
    \ENDIF
\ENDFOR

\STATE $\mathbf{x} \gets \mathbf{x}_{k+1}$
\end{algorithmic}
\end{algorithm}

\subsection{Criterion for switching precision for $\Br_k$ and $\Bq_k$}
\label{sec: criterion rq}

To achieve the target tolerance $\tol$,  the inequality \eqref{eq:conv_true_res} must be satisfied.
% \begin{equation}
% \|\Bzeta_k\|   + \|\Br_k\| 
%    \leq  \tol\,\|\Bb\|
%   \label{eq: tol bound zeta}
% \end{equation}
Therefore, we aim to derive a bound on $\|\Bzeta_k\|$ \eqref{eq:res_gap} to ensures the prescribed accuracy is met. 
We begin by establishing the following results under three assumptions:
\begin{assumption}
There exists a constant \( C > 0 \), independent of the matrix size \( n \), such that for all \( k \),
\begin{equation}
    \left\| \mathrm{fl}(\BA \Bp_k, \epsqk) - \BA \Bp_k \right\| \leq C \epsqk \left\| \BA \Bp_k \right\|.
    \label{eq:bound Ap}
\end{equation}
\label{asp:1}
\end{assumption}
\begin{assumption}
     There exists a constant $C_x>0$ such that $\|\Bx_k\| \leq C_x \|\Bx\|$ for all $k$.
     \label{asp:3}
\end{assumption}
\begin{assumption}
    Let $\bar{k}$ denote the total number of iterations required by the algorithm. Then $\bar{k} \epsxk \|\BA\| \|\Bx\|$ is sufficiently small.
\label{asp:2}
\end{assumption}
\begin{remark}
\Cref{asp:1} adopts a less conservative error bound than is commonly used in the literature:
\[
\left\| \mathrm{fl}(\BA \Bp_k, \epsqk) - \BA \Bp_k \right\| \leq \hat{C} \, \epsqk \|\BA\| \|\Bp_k\|,
\]
where \( \hat{C} = m_r n^{1/2} \) in \cite{golub2013matrix} and \( m_r \) denotes the maximum number of nonzeros per row. While sharper bounds can be obtained through probabilistic rounding error analysis \cite{higham2002accuracy, higham2019new}, we adopt \eqref{eq:bound Ap} for practicality, as the more rigorous alternatives often yield  pessimistic estimates. 
In practice, $C=1$ seems to yield satisfactory results, although it may be somewhat optimistic.

\Cref{asp:3} holds when the initial guess  \( \Bx_0 \) is comparable in the size to the exact solution \( \Bx \), which is often the case in practice \cite{hestenes1952methods, greenbaum1997estimating}.

\Cref{asp:2} reflects the fact that \( \Bx_k \) is computed in double precision, and hence \( \epsxk \) corresponds to the unit roundoff for FP64. Typically, both \( \|\BA\| \) and the total iteration count \( \bar{k} \) are moderate, ensuring that the term \( \bar{k} \, \epsxk \|\BA\| \|\Bx\| \) is negligible in comparison with the numerical errors arising from the residual updates.
\end{remark}

We now derive a bound for \( \|\Btheta_k\| \), which is defined in \cref{eq:res perturb}. Since both \( \Br_k \) and \( \Bq_k \) are computed in precision \( u_{r,k} \), it follows that \( \epsqk = \epsrk \) in \cref{eq:pertub bound}. From \cref{asp:1}, we obtain
\begin{equation}
\begin{aligned}
    \|\Btheta_k\| & \leq   \varepsilon_{r,k} \|\Br_{k-1}\| +
    (2\varepsilon_{r,k} + \varepsilon_{r,k}^2)\|\alpha_{k-1} \BA \Bp_{k-1}\| +
    (1+\varepsilon_{r,k})^2 \left\| \alpha_{k-1} \left( \mathrm{fl}(\BA \Bp_{k-1}, \epsrk) - \BA \Bp_{k-1} \right) \right\| \\
    &\leq \varepsilon_{r,k} \|\Br_{k-1}\| +
    (2\varepsilon_{r,k} + \varepsilon_{r,k}^2)\|\alpha_{k-1} \BA \Bp_{k-1}\| +
    C (1+\varepsilon_{r,k})^2 \varepsilon_{r,k} \|\alpha_{k-1} \BA \Bp_{k-1}\| \\
    &\leq \varepsilon_{r,k} \|\Br_{k-1}\| +
    \left( 2\varepsilon_{r,k} + \varepsilon_{r,k}^2 + C\, \varepsilon_{r,k} (1+\varepsilon_{r,k})^2 \right) \|\alpha_{k-1} \BA \Bp_{k-1}\| \\
    % &\leq \varepsilon_{r,k} \|\Br_{k-1}\| +
    % \left( \lp 2+C \rp\varepsilon_{r,k} + \lp 1+3C \rp\varepsilon_{r,k}^2 \right) \|\alpha_{k-1} \BA \Bp_{k-1}\| \\
    &\leq \varepsilon_{r,k} \|\Br_{k-1}\| +
    \left( \lp 2+C \rp\varepsilon_{r,k} + \lp 1+3C \rp\varepsilon_{r,k}^2 \right) 
    \left( \|\Br_{k-1}\| + \|\Br_k\| + \| \Btheta_k\| \right).
\end{aligned}
\end{equation}
As long as $ 1 -    \left( \lp 2+C \rp\varepsilon_{r,k} + \lp 1+3C \rp\varepsilon_{r,k}^2 \right) >0 $, this can be written in the form
\begin{equation}
\begin{aligned}
\|\boldsymbol{\theta}_k\| 
& \leq \lp 3+ C \rp \varepsilon_{r,k} \|\Br_{k-1}\| +  \lp 2+ C \rp \varepsilon_{r,k} \|\Br_k\| + \mathcal{O}(\varepsilon_{r,k}^2) \left( \|\mathbf{r}_{k-1}\| + \|\mathbf{r}_k\| \right) .
\end{aligned}
\label{eq:bound theta}
\end{equation}
Similarly, with \cref{asp:3}, 
$\|\Bxi_k\|$ can be bounded  as follows:
\begin{equation}
\begin{aligned}
\left\|\Bxi_k\right\|
 \leq \epsxk \left(3\left\|\Bx_{k-1}\right\|+2\left\|\Bx_{k}\right\|\right)+O\left(\epsxk^{2}\right)\left(\left\|\Bx_{k-1}\right\|+\left\|\Bx_k\right\|\right) \leq 5 \,k \,\epsxk \|\Bx\| + O\left(\epsxk^{2}\right) \|\Bx\|.
\end{aligned}
\label{eq:bound xi}
\end{equation}
From \eqref{eq:res_gap}, \eqref{eq:bound theta}, and \eqref{eq:bound xi}, we obtain the following bound for  $\|\Bzeta_k\|$: 
\begin{equation}
\begin{aligned}
    \|\Bzeta_k\| 
    &\leq \|\Bb - \BA\Bx_0 - \Br_0\| + \|\BA\| \sum\limits_{t=1}^{k} \|\Bxi_t\| + \sum\limits_{t=1}^{k}  \|\Btheta_t\| \\
    & \leq  5\, k\, \epsxk \|\BA\|\|\Bx \| + 
    \sum\limits_{t=0}^{k} \epsrk \lp \lp 3+C \rp \|\Br_{t-1}\| + \lp 2+C \rp \|\Br_t\|  \rp + O\left(\varepsilon_{x,t}^{2}\right) + O\left(\varepsilon_{r,t}^{2}\right) \\
    & \approx  \sum\limits_{t=0}^{k} \epsrk \lp \lp 3+C \rp \|\Br_{t-1}\| + \lp 2+C \rp \|\Br_t\|  \rp.
\end{aligned}
\label{eq: bound zeta}
\end{equation}
The last approximation follows by neglecting the first term  (\cref{asp:2}) and  the second order terms.
The norm $\|\Bzeta_k\|$ estimates the best attainable accuracy for sufficiently large $k$  (so that $\|\Br_k\|$ is small), as shown in \eqref{eq:conv_true_res}.
Switching precision $u_{r,k}$ from FP64 to FP32 at iteration ${k}$, we define an indicator $\widehat{\eta}_m( {k} )$ to estimate the best attainable accuracy at iteration $m$:
\begin{equation}
   \widehat{\eta}_m( {k} ) :=   \sum\limits_{t=0}^{ {k}} \epsD \lp \lp 3+C \rp \|\Br_{t-1}\| + \lp 2+C \rp \|\Br_t\|  \rp + \sum\limits_{t={k}}^{ m} \epsS \lp \lp 3+C \rp \|\Br_{t-1}\| + \lp 2+C \rp \|\Br_t\|  \rp,
    \label{eq: Weta_k}
\end{equation}
where  $\epsS$ and $\epsD$ are the unit roundoff errors of single and double precision, respectively, and  $m$ denotes the final iteration index.
In practice, however, $\widehat{\eta}_m( {k} )$ cannot be evaluated until the run completes, since it require all residual norms $\|\Br_t\|$ for ${k}<t\leq m$.
Therefore, we introduce a computable heuristic indicator that predicts the attainable accuracy without relying on residual norms $\|\Br_t\|$ beyond iteration $k$.
Let integer $d$ be a fixed delay parameter. 
If reducing precision at iteration $k-d$ satisfies the condition
 $\widehat{\eta}_m( k-d )\leq \tol \, \|\Bb\|$, then switching at $k$ will also satisfy the tolerance (since $\widehat{\eta}_m(k-d) > \widehat{\eta}_m(k)$).
Moreover, we assume the residual decays rapidly after these $d$ iterations, such that 
\begin{equation}
    \sum\limits_{t=k+1}^{m} \|\Br_{{t}}\| \ll \sum\limits_{t=k-d}^{k} \|\Br_{{t}}\|.
    \label{eq: res norm decay}
\end{equation}
Now, we define the computable  attainable accuracy indicator as the dominant contribution in  $\widehat{\eta}_m(k-d)$, expressed as: 
\begin{equation}
\eta_k := 
\sum\limits_{t=k-d}^{k} \epsS \left( \lp 3+C \rp\|\Br_{t-1}\| + \lp2 + C\rp \|\Br_t\| \right).
\label{eq: eta_k}
\end{equation}
Additionally, we  define the criterion for switching the precision of $\Br_{k}$ and $\Bq_{k}$ from {FP64} to {FP32}:
\begin{equation}
    \eta_k \leq {\tol}\, \|\Bb\|\,.
    \label{eq: precision switch criterion}
\end{equation}
In our numerical experiments (see \cref{sec: precision of r and q}), we demonstrate that choosing $d=10$ and $C=1$ in \cref{eq: eta_k} provide a robust and effective compromise. This delay parameter is inspired by a similar approach used for estimating the $\BA$-norm of the error in \cite{hestenes1952methods}.

When the convergence rate is nearly linear, the 
residual norms can be estimated without
a delay parameter.
We can estimate the average rate of convergence 
by 
$\rho = (\|\Br_k\| / \|\Br_{k-\ell}\|)^{1/\ell}$, for a chosen parameter $\ell$, and estimate
$\widehat{\eta}_{m}(k)$ using 
$\|\Br_t\| = \|\Br_{k}\|\rho^{t-k}$.
Neglecting all terms proportional to $\epsD$, we therefore define our attainable accuracy indicator $\eta_k$ and the corresponding precision switching criterion as follows:
\begin{equation}
\begin{aligned}
  \eta_k & = 
  \sum\limits_{t={{k}}}^{m} \epsS \left( \lp 3+C\rp\|\Br_{t-1}\| + \lp2+C\rp \|\Br_t\| \right)\\
  & \leq 
\epsS \lp 5 + 2C\rp \| \Br_{k-1}\|  \frac{1}{1-\rho}   \leq  {tol} \, \|\Bb\|.
\end{aligned}    
\label{eq:linear criterion}
\end{equation}
Therefore, when the convergence is almost linear, we use \cref{eq:linear criterion} to choose the precision switching step instead of \cref{eq: eta_k}. This criterion is tested in \cref{sec: well conditioned matrices} and \cref{sec: performance}.

Substituting \( \eta_k \) as an estimate for \( \max_{t\leq m} \|\Bzeta_t\| \) in \eqref{eq:conv_true_res} and combining this with the stopping criterion (line 18 in \cref{alg:MPDS-PCG}), we obtain a bound for the true residual norm:
\[
\|\Bb - \BA\Bx_m\| \lesssim \eta_k + \|\Br_m\| \leq 2\, \tol\,\|\Bb\|.
\]
Thus, given a target accuracy \( \tol \), we can choose when to switch the precision based on either \cref{eq: eta_k} or \cref{eq:linear criterion}, such that the true residual norm satisfies the accuracy.

\section{Numerical Results}
\label{sec:numerical}

In this section, we test the performance of   \APCG{} (Algorithm~\ref{alg:MPDS-PCG}) on the solution of a set of example linear systems.
We begin by revisiting the motivating example in \cref{sec: revisit motivating}. Then in \cref{sec: precision of r and q}, we assess the attainable accuracy of \APCG{}, and then in \cref{sec: well conditioned matrices} and \cref{sec: matcol} we examine its convergence behavior.
We evaluate the performance of \APCG{} within a large-scale GPU-accelerated high-order finite element solver in \cref{sec: performance}. 
Finally, we illustrate some limitations of the precision selector \cref{eq:precision_selector} in \cref{sec: limitation} through several  experiments. A summary of all examples discussed in this section is provided in \cref{tab:all examples}.
\begin{table}[htpb]
    \centering
    \renewcommand{\arraystretch}{1.3}
    \setlength{\tabcolsep}{8pt}
    \begin{tabular}{@{} l l  r c c @{}}
    \toprule
    {Example} & Description  & {\(n\)} & {Est. \(\kappa(\BM^{-1}\BA)\)} & {Convergence} \\
    \midrule
    \ref{ex: motivating} & Motivating example  & 100 & \(10^3\) & Superlinear \\
    \ref{ex:bcsstk}&  \texttt{bcsstk08.mat}  & 1074 & \(3.77\times10^3\) & Superlinear \\
    \ref{ex:622bus}&  \texttt{622\_bus.mat}  & 662 & \(4.46\times10^4\) & Superlinear \\
    \ref{ex: linear conv1}&  Evenly distributed eigen.  & 101 & 10 & Linear \\
    \ref{ex: linear conv2}&  FEM with multigrid prec.  & \(9.1\times10^6\) & 187.2 & Linear \\
    \ref{ex: not large outlying eigs}&  Evenly distributed eigen. & 101 & \(10^3\) & Superlinear \\
    \ref{ex: jac libp}&   FEM with Jacobi prec.  & \(2.2\times10^7\) & $1.99\times 10^3$ & Linear \\
    \ref{ex: limit res}&  \texttt{nos6.mat}  & 675 & \(3.49\times10^6\) & Superlinear \\
    \ref{ex: limitation 1}&  Illustrative example  & 10 & \(10^3\) & Oscillatory \\
    \ref{ex: limitation 2}&  Many large outlying eigen. & 50 & \(10^3\) & Oscillatory \\
    \bottomrule
    \end{tabular}
    \caption{Summary of test examples. \(n\) is the matrix size, \(\kappa\) is the estimated condition number of the preconditioned matrix, and the last column indicates observed convergence behavior.}
    \label{tab:all examples}
\end{table}

In this section, unless otherwise specified, in the precision selector function \cref{eq:precision_selector}, $\eta_k$ 
  is defined by \cref{eq: precision switch criterion}  with $d=10$ and $C=1$ under superlinear convergence and by \cref{eq:linear criterion} with $\ell=5$ 
  under linear convergence.
  The following heuristic thresholds are used: 
\begin{equation}
    \tau_{z,s} = 10^{-4}, \quad \tau_{z,h} = 10^{-6}.
    \label{eq: criterion}
\end{equation}
In the following subsections, we compare several PCG variants:
\begin{enumerate}
    \item \textbf{PCG(FP64)}: All vectors are computed and stored in double precision.
    \item \textbf{\APCG(FP64)}: Initially uses FP64 for both \( u_{z,k} \) and \( u_{r,k} \), which are then adaptively reduced during the iterations.
    \item \textbf{\APCG(FP16)}: Use \( \Bz_k \) and \( \Bp_k \)  in FP16, while \( u_{r,k} \) is adaptively reduced.
    \item \textbf{\APCG($\urk$) }: Keeps \( u_{z,k} \) fixed in FP64, and adaptively reduces \( u_{r,k} \); used to assess the effectiveness of the precision-switching criterion for \( u_{r,k} \).
    \label{enu:pcg_variant}
\end{enumerate}

Moreover, unless stated otherwise, the right hand side vector $\Bb$ is randomly generated,  the initial guess $\Bx_0$ is the zero vector, the preconditioner is the identity matrix, 
and the stopping tolerance $\tol = 10^{-10}$.
In the numerical results, the updated residual corresponds to the residual computed within the \APCG{}  iteration (line~16 in Algorithm~\ref{alg:MPDS-PCG}), whereas the true residual is defined as $ \Bb - \BA \Bx_k$.

\subsection{Motivating example}
\label{sec: revisit motivating}
As demonstrated in the motivating example (see \cref{ex: motivating} and \cref{fig: motivation example}), applying the PCG algorithm in single or half precision leads to reduced attainable accuracy and slower convergence relative to PCG in double precision; half precision in particular also suffers from numerical underflow.

\begin{figure}[htbp]
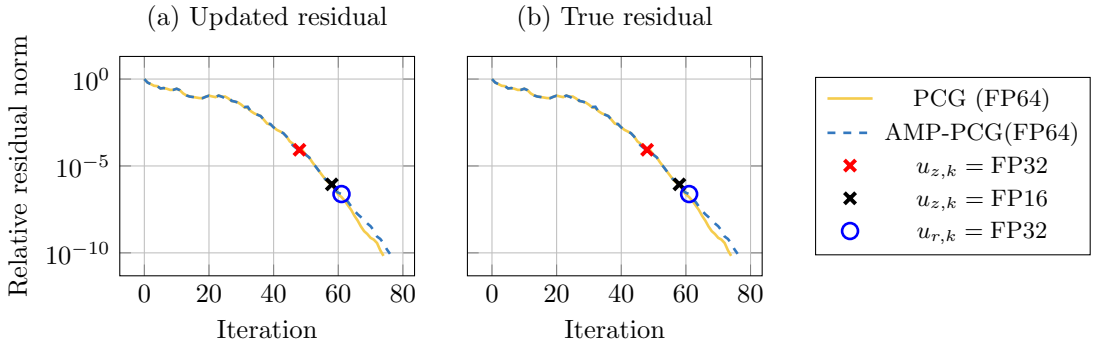

\centering
\begin{tikzpicture}
\begin{groupplot}[
    group style={
        group name=myplot,
        group size=2 by 1,
        horizontal sep=20pt,
        vertical sep=2cm
    },
    height=4.5cm,
    width=5.5cm
]
  % First plot with y-axis labels
  \nextgroupplot[
  title={(a) Updated residual},
    xlabel={Iteration},
    ylabel={Relative residual norm},
    ymode=log,
    % ymin=1e-10,
    ymax=20,
        ytick={1e0, 1e-5, 1e-10, 1e-14},
     grid=both,
  ]
  \input{5.Numerical/data/motivating_after_1}
  
  % Second plot without y-axis tick labels
  \nextgroupplot[
    title={(b) True residual},
    xlabel={Iteration},
    ymode=log,
    % ymin=1e-9,
    ymax=20,
    yticklabels={},
    ytick={1e0, 1e-5, 1e-10, 1e-14},
     grid=both,
  ]
  \input{5.Numerical/data/motivating_after_2}
 
\end{groupplot}

\path (myplot c1r1.south west|-current bounding box.south)--
      coordinate(legendpos)
      (myplot c2r1.south east|-current bounding box.south);

\matrix[
    matrix of nodes,
    anchor=west, % align the left edge of the legend with the anchor point
    draw,
    inner sep=0.2em,
    nodes={font=\small}
] at ([xshift=2em]myplot c2r1.east) % shift a bit to the right
{
    \ref{plot:fp64} & PCG (FP64) \\
    \ref{plot:fp64mixed} & \APCG(FP64) \\
    \ref{plot:fp64tofp32} & $ u_{z,k}=\mathrm{FP32}$  \\
    \ref{plot:fp32tofp16} & $ u_{z,k}=\mathrm{FP16}$ \\
        \ref{plot:residualfp64tofp32} & $ u_{r,k}=\mathrm{FP32} $
    \\
};
\end{tikzpicture}
\caption{
Revisiting the motivating example (see \cref{ex: motivating,fig: motivation example}): 
convergence of (a) the relative updated residual norm and (b) the relative true residual norm for PCG(FP64) and \APCG(FP64), as discussed in \cref{sec: revisit motivating}. Cross markers indicate the iterations where the precision \( u_{z,k} \), used for \( \mathbf{z}_k \) and \( \mathbf{p}_k \), is reduced to single and then to half precision. The blue circle marks where \( u_{r,k} \), used for \( \mathbf{r}_k \) and \( \mathbf{q}_k \), is reduced from double to single precision.
}

\label{fig: motivation example again}
\end{figure}

We now apply  our \APCG{} algorithm to the same example. 
In \APCG{}(FP64), all vectors $\mathbf{z}_k$, $\mathbf{p}_k$, $\mathbf{r}_k$, and $\mathbf{q}_k$ are initially computed and stored in double precision. Their precisions are adaptively reduced based on the precision selector function \cref{eq:precision_selector}. In \cref{fig: motivation example again}, the red and black cross markers  denote when the precision $u_{z,k}$, used for $\mathbf{z}_k$ and $\mathbf{p}_k$, is reduced first to single precision and then to half precision, respectively; the blue circle markers denote when precision $u_{r,k}$, associated with $\mathbf{r}_k$ and $\mathbf{q}_k$, is reduced from double to single precision.

The convergence of both the updated and true residuals under \APCG{}(FP64) closely mirrors that of the  PCG(FP64). This demonstrates the effectiveness of the proposed \APCG{} in maintaining convergence and accuracy while using low precision. In the next subsections, we further evaluate the performance of \APCG{} from multiple perspectives.

\subsection{Impact of the precision of $\Br_k$ and $\Bq_k$ on attainable accuracy}
\label{sec: precision of r and q}
We assess whether the mixed precision strategy for $\Br_k$ and $\Bq_k$ achieves a final residual norm within the prescribed tolerance \(\tol\).
All other vectors remain in double precision. We test \cref{ex: motivating} and two SPD matrices 
from the SuiteSparse Matrix Collection \cite{davis2011university}, \texttt{bcsstk08.mat} and \texttt{622\_bus.mat}, which are commonly used for testing PCG. These two matrices are scaled on both sides by the square root of the inverse of its diagonal entries.  
\begin{example}\label{ex:bcsstk}
Consider the $1074\times1074$ matrix \texttt{bcsstk08.mat}, for which  the condition number after scaling is approximately  $3.77\times10^3$.
\end{example}
\begin{example}\label{ex:622bus}
Consider the $662\times662$ matrix \texttt{622\_bus.mat}, for which  the condition number after scaling is about $4.46\times10^4$.
\end{example}
We consider two options for the stopping tolerance $\tol$  used in line 18 of \cref{alg:MPDS-PCG}:
\[
\tol = 10^{-6}\, \text{ and } \, \tol=10^{-8}.
\]

\begin{figure}[h!]
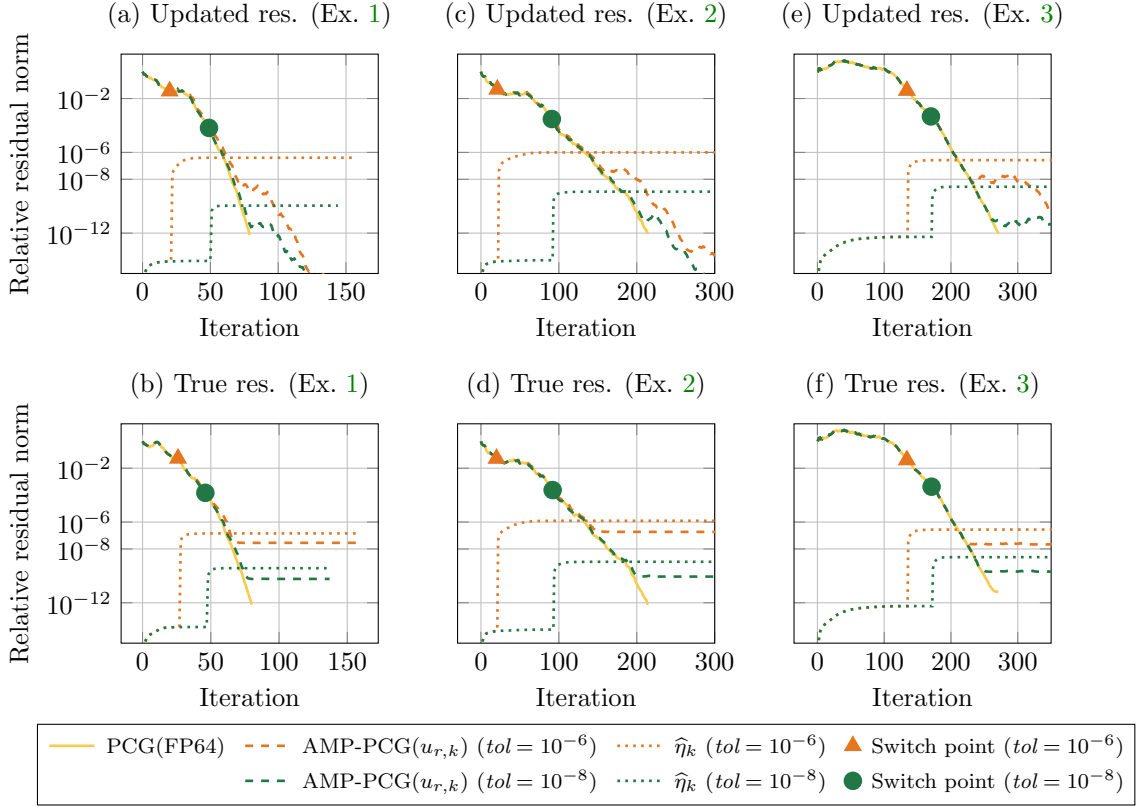

\centering
\begin{tikzpicture}
\begin{groupplot}[
    group style={
        group name=myplot,
        group size=3 by 2,
        horizontal sep=30pt,
        vertical sep=2cm
    },
    height=4.5cm,
    width=5cm,
]
  % First plot with y-axis labels
  \nextgroupplot[
  title={(a) Updated res. (\exref{ex: motivating})},
    xlabel={Iteration},
    ylabel={Relative residual norm},
    ymode=log,
    ymin=1e-15,
    ymax=20,
    grid=both,
     ytick={1e-2, 1e-6, 1e-8,1e-12},
  ]
  \input{5.Numerical/data/resPrec_1_1}

  % First plot with y-axis labels
  \nextgroupplot[
  title={(c) Updated res. (\exref{ex:bcsstk})},
    xlabel={Iteration},
    ymode=log,
    ymin=1e-15,
    ymax=20,
    xmax=300,
    grid=both,
     ytick={1e-2, 1e-6, 1e-8,1e-12},
     yticklabels={},
  ]
  \input{5.Numerical/data/resPrec_2_1}

   % First plot with y-axis labels
  \nextgroupplot[
  title={(e) Updated res. (\exref{ex:622bus})},
    xlabel={Iteration},
    ymode=log,
    ymin=1e-15,
    ymax=20,
        xmax=350,
    grid=both,
     ytick={1e-2, 1e-6, 1e-8,1e-12},
    yticklabels={},
  ]
  \input{5.Numerical/data/resPrec_3_1}
  
  % Second plot without y-axis tick labels
  \nextgroupplot[
    title={(b) True res. (\exref{ex: motivating})},
    xlabel={Iteration},
    ylabel={Relative residual norm},
    ymode=log,
    ymin=1e-15,
    ymax=20,
    grid=both,
     ytick={1e-2, 1e-6, 1e-8,1e-12},
  ]
  \input{5.Numerical/data/resPrec_1_2}

  % Second plot without y-axis tick labels
  \nextgroupplot[
    title={(d) True res. (\exref{ex:bcsstk})},
    xlabel={Iteration},
    ymode=log,
    ymin=1e-15,
    ymax=20,
        xmax=300,
    grid=both,
     ytick={1e-2, 1e-6, 1e-8,1e-12},
   yticklabels={},
  ]
  \input{5.Numerical/data/resPrec_2_2}

  % Second plot without y-axis tick labels
  \nextgroupplot[
    title={(f) True res. (\exref{ex:622bus})},
    xlabel={Iteration},
    ymode=log,
    ymin=1e-15,
    ymax=20,
        xmax=350,
    grid=both,
     ytick={1e-2, 1e-6, 1e-8,1e-12},
     yticklabels={},
  ]
  \input{5.Numerical/data/resPrec_3_2}
  
\end{groupplot}

\path (myplot c1r1.south west|-current bounding box.south)--
      coordinate(legendpos)
      (myplot c3r1.south east|-current bounding box.south);

\matrix[
    matrix of nodes,
    anchor=south,
    draw,
    inner sep=0.2em,
    nodes={font=\footnotesize}
] at ([yshift=-8ex]legendpos)
{
    \ref{plot:fp64} & PCG(FP64) & [3pt] & 
    \ref{plot:resPrec1e6} & \APCG($\urk$) $(\tol = 10^{-6})$ & [3pt] &
    \ref{plot:resPrec1e6_ind} & $\widehat{\eta}_k$ $(\tol = 10^{-6})$  & [3pt] &
    \ref{plot:residual2} & Switch point $(\tol = 10^{-6})$  \\
    & & [3pt] & \ref{plot:resPrec1e8} & \APCG($\urk$) $(\tol = 10^{-8})$ & [3pt] &
    \ref{plot:resPrec1e8_ind} & $\widehat{\eta}_k$ $(\tol = 10^{-8})$ & [3pt] &
    \ref{plot:residual3} & Switch point $(\tol = 10^{-8})$ \\
};
\end{tikzpicture}
\caption{
Convergence curves of the relative updated residual norm ((a), (c), and (e)) and the relative true residual norm ((b), (d), and (f)) for PCG(FP64) and \APCG($\urk$)   applied to \cref{ex: motivating,ex:bcsstk,ex:622bus}, respectively. Two stopping tolerances, \(10^{-6}\) and \(10^{-8}\), are tested. 
At the switch point, the precision of \( \Br_k \) and \( \Bq_k \) is reduced from FP64 to FP32; all other vectors remain in FP64. The dotted curves show the indicator \( \widehat{\eta}_k \) (see \cref{eq: Weta_k}), which provides an estimate of the final attainable accuracy. 
}
\label{fig:vary taur}
\end{figure}

According to the precision selector function \cref{eq:precision_selector}, varying \(\tol\) changes when the precision  \(u_{r,k}\) is reduced, thereby impacting the attainable accuracy.
In \cref{fig:vary taur}, we present convergence curves for both the updated residual (top row) and the true residual (bottom row) for PCG and \APCG($\urk$) . In these experiments, the iterations are not terminated based on the stopping criterion in \cref{alg:MPDS-PCG}, allowing us to examine the final attainable residual accuracy. Markers indicate the iteration at which the precision switch occurs for each value of $\tol$.
The results indicate that using low precision for $\Br_k$ and $\Bq_k$ does not affect the convergence rate in these examples; however, it does limit the attainable accuracy. 
As illustrated in subfigures (b), (d), and (f), 
the precision selector function successfully drives the true residual below the prescribed tolerance  $\|\Br_k\| \leq \tol\, \|\Bb\|$.
Additionally, we compute the indicator $\widehat{\eta}_k$, defined in \cref{eq: Weta_k}.
This indicator is an estimate of the final attainable accuracy, and its values are shown as dotted lines. 
In these examples, $\widehat{\eta}_k$ closely approximates the attainable accuracy, confirming the effectiveness of the indicator.

The results show that using the precision selector function \cref{eq:precision_selector} to determine when to reduce the precision of \( \Br_k \) and \( \Bq_k \) is effective in achieving the target accuracy.
Since the updated and true residual norms closely agree prior to satisfying the stopping criterion, we report only the updated residual norm in subsequent subsections.

\subsection{Cases with linear convergence}
\label{sec: well conditioned matrices}
We consider two matrices that either have a small condition number or are equipped with an effective preconditioner, such that the condition number of the preconditioned system is small. In these cases, we anticipate linear convergence. As discussed in \cref{sec: z and p precision}, within the linear convergence regime, using \( \Bp_k \) and \( \Bz_k \) in half precision has negligible impact on the convergence rate. Therefore, we fix \( u_{z,k} = \mathrm{FP16} \), using \( \Bz_k \) and \( \Bp_k \) in half precision throughout the iterations.
\begin{example}
\label{ex: linear conv1}
Define  $\BA \in \mathbb{R}^{(n+1)\times (n+1)}$ with uniformly distributed eigenvalues
\[
\lambda_i = \frac{i}{10n} + \frac{n-i}{n}, \quad i=0,\dots,n,
\]
where $n=100$. The eigenvector matrix is a random orthonormal matrix. 
\label{ex:small cond 1}
\end{example}
\begin{example}
\label{ex: linear conv2}
We consider a matrix arising from a high-order finite element discretization of the Poisson equation:
\begin{equation}
- \Delta u(x,y,z) = f(x,y,z) \quad \text{in } \Omega=[-0.5,0.5]^3,
\label{eq:poisson}
\end{equation}
subject to homogeneous Dirichlet boundary conditions. The forcing function is given by
\begin{equation}
    f(x,y,z)= 1 + 3\pi^2 \sin(\pi x)\sin(\pi y)\sin(\pi z).
\label{eq:f}
\end{equation}
We discretize the equation in the weak form using high-order finite element method with 
a polynomial degree of $N=7$ on a  $30\times30\times30$ hexahedral  Kershaw mesh with $\varepsilon_{\text{Kershaw}}=0.3$, which is used as the mesh of benchmark problems by the Center for Efficient Exascale Discretization  \cite{kolev_2021_4672664}. 
The total degrees of freedom (DoFs) is about $9\times 10^6$.
We solve the system using both PCG and \APCG(FP16), preconditioned with a \( p \)-multigrid with a second-order Chebyshev-accelerated Jacobi smoother \cite{karakus2019gpu}. On the coarsest level (\( N = 1 \)), an algebraic multigrid solver is applied. 
The condition number of the preconditioned matrix is about 187, estimated from Ritz value of iteration 153. 
All experiments are performed using \textit{libParanumal} \cite{ChalmersKarakusAustinSwirydowiczWarburton2020} on a NVIDIA GTX 4090 GPU. 
\end{example}

% \begin{figure}[htbp]
%     \centering
%     \includegraphics[width=0.3\linewidth]{5.Numerical/fig/mesheps03.jpeg}
%     \caption{Kershaw mesh used in \cref{ex: linear conv2}.}
%     \label{fig:ex2mesh}
% \end{figure}

\begin{figure}[htpb]
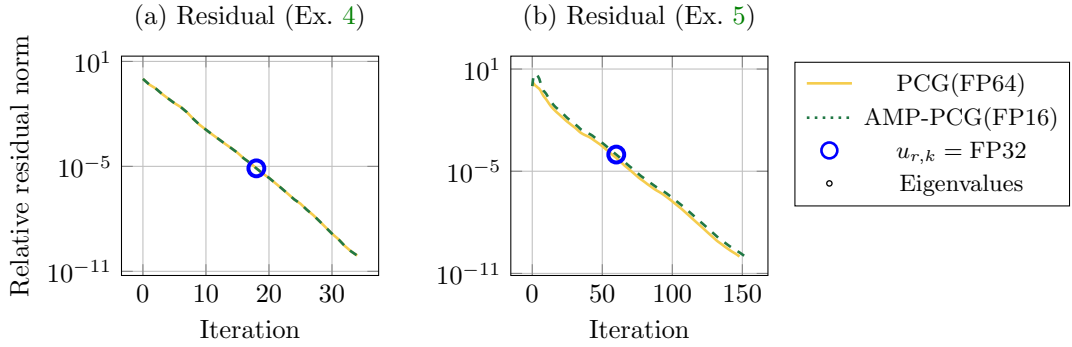

\centering
\begin{tikzpicture}
\begin{groupplot}[
    group style={
        group name=myplot,
        group size=2 by 1,
        horizontal sep=50pt,
        vertical sep=2cm
    },
    height=4.5cm,
    width=5cm
]
  % (a) Residual - Example 1
  \nextgroupplot[
    title={(a) Residual (\exref{ex: linear conv1})},
    xlabel={Iteration},
    ylabel={Relative residual norm},
    ymode=log,
    ymax=20,
    grid=both,
  ]
  \input{5.Numerical/data/smallCond_1}

  % (b) Eigenvalues - Example 1
  \nextgroupplot[
    title={(b) Residual (\exref{ex: linear conv2})},
    xlabel={Iteration},
        ymode=log,
        grid=both,
  ]
  \input{5.Numerical/data/smallCond_2}

  % % (c) Residual - Example 2
  % \nextgroupplot[
  %   title={(c) Residual (\exref{ex: linear conv2})},
  %   xlabel={Iteration},
  %   ylabel={Relative residual norm},
  %   ymode=log,
  %   grid=both,
  % ]
  % \input{5.Numerical/data/smallCond_2}

  % % (d) Ritz values - Example 2
  % \nextgroupplot[
  %   title={(d) Ritz values (\exref{ex: linear conv2})},
  %   xlabel={Index},
  %   ylabel={Ritz values},
  %   ymode=log,
  % ]
  % \input{5.Numerical/data/smallCond_2eig}
\end{groupplot}

% Legend positioning
\matrix[
    matrix of nodes,
    anchor=west,
    draw,
    inner sep=0.2em,
    nodes={font=\small}
] at ([xshift=1em, yshift=1.2em]myplot c2r1.east) {
    \ref{plot:fp64} & PCG(FP64) \\
    \ref{plot:fp16} & \APCG(FP16) \\
    \ref{plot:residualfp64tofp32} & $u_{r,k} = \mathrm{FP32}$ \\
    \ref{plot:eig} & Eigenvalues \\
};

\end{tikzpicture}
\caption{
Convergence of the relative residual norm for PCG(FP64) and \APCG(FP16) on (a) \cref{ex: linear conv1}  and (b) \cref{ex: linear conv2}.
In \APCG(FP16),  \( \Bz_k \) and \( \Bp_k \) are in half precision throughout the iterations.
The blue circle marks where \( u_{r,k} \), used for \( \mathbf{r}_k \) and \( \mathbf{q}_k \), is reduced from double to single precision.
}
\label{fig:small cond}
\end{figure}

In \cref{fig:small cond}, we present the convergence curves for PCG(FP64) and \APCG(FP16), where \( \Bp_k \) and \( \Bz_k \) are stored in FP16 throughout all iterations, while the precision of \( \Br_k \) and \( \Bq_k \) is adaptively reduced according to \cref{eq:precision_selector} with $\eta_k$ defined in \cref{eq:linear criterion}.
Subfigure~(a) and (b), corresponding to the precision selector \cref{ex: linear conv1} and  \cref{ex: linear conv2}, respectively, shows that the convergence of \APCG(FP16) is nearly identical to that of  PCG(FP64). The blue circle indicates the iteration at which the precision of \( \Br_k \) and \( \Bq_k \) is reduced to single precision.

These results demonstrate that when the convergence is linear, using half precision for \( \Bz_k \) and \( \Bp_k \) rarely affect the convergence rate, which is consistent with the analysis presented in \cref{sec: z and p precision}.

\subsection{Matrices without large, outlying eigenvalues}
\label{sec: matcol}
We consider matrices whose eigenvalue distributions do not contain large, outlying eigenvalues. These matrices typically exhibit linear convergence after a few initial iterations, allowing the precision of \( \Bz_k \) and \( \Bp_k \) to be  reduced once the linear regime is reached. 
In contrast, matrices with many large, outlying eigenvalues often lead to oscillatory convergence in the residual norm, as demonstrated in \cref{sec: choice of tauz}. A similar distinction between matrices with and without large, outlying eigenvalues is discussed in \cite{greenbaum2021convergence, carson2024towards}.
In this subsection, we examine \cref{ex:bcsstk}, \cref{ex:622bus}, and a synthetic matrix with uniformly distributed eigenvalues, defined as follows.
\begin{example}
We consider a matrix $\BA \in \mathbb{R}^{(n+1)\times (n+1)}$ with uniformly distributed eigenvalues
\[
\lambda_i = \frac{i}{1000n} + \frac{n-i}{n}, \quad i=0,\dots,n,
\]
with $n=100$ from $\lambda_n = 10^{-3}$ to $\lambda_0 = 1$. The eigenvector matrix is a random orthonormal matrix. 
\label{ex: not large outlying eigs}
\end{example} 

\begin{figure}[h!]
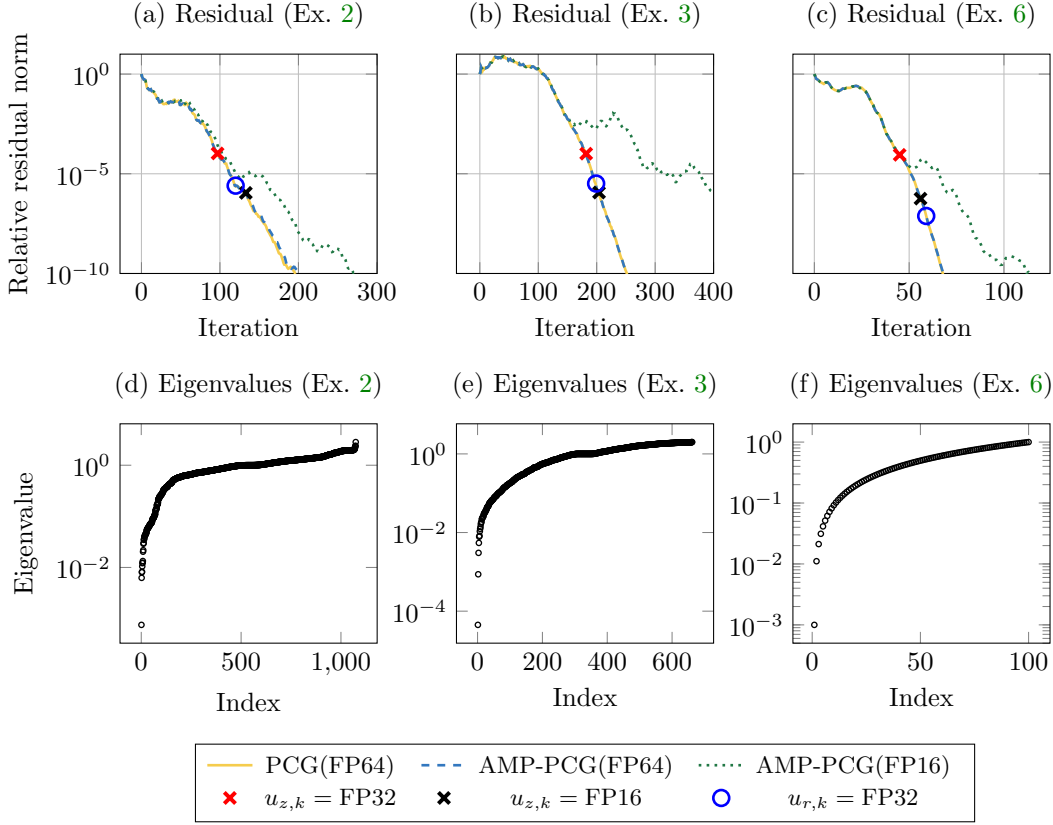

\centering
\begin{tikzpicture}
\begin{groupplot}[
    group style={
        group name=myplot,
        group size= 3 by 2,
        horizontal sep=30pt,
        vertical sep=2cm
    },
    height=4.5cm,
    width=5cm
]

  % 4 plot with y-axis tick labels
  \nextgroupplot[
    title={(a)  Residual (\exref{ex:bcsstk})},
    xlabel={Iteration},
    ylabel={Relative residual norm},
    ymode=log,
    ymin=1e-10,
    ymax=10,
    grid=both,
  ]
  \input{5.Numerical/data/okCond_2_1}

    % 7 plot with y-axis tick labels
  \nextgroupplot[
    title={(b) Residual (\exref{ex:622bus})},
    xlabel={Iteration},
    ymode=log,
    ymin=1e-10,
    ymax=10,
    xmax=400,
    yticklabels={},
    grid=both,
  ]
  \input{5.Numerical/data/okCond_3_1}

  % First plot with y-axis labels
  \nextgroupplot[
    title={(c) Residual (\exref{ex: not large outlying eigs})},
    xlabel={Iteration},
    ymode=log,
    ymin=1e-10,
    ymax=10,
    grid=both,
        yticklabels={},
  ]
  \input{5.Numerical/data/okCond_1_1}
  
    % 6 plot without y-axis tick labels
  \nextgroupplot[
    title={(d) Eigenvalues (\exref{ex:bcsstk})},
    xlabel={Index},
    ylabel={Eigenvalue},
    ymode=log,
    % ymin=1e-10,
    % ymax=10,
    % yticklabels={}
  ]
  \input{5.Numerical/data/okCond_2_eig}
  
    % 9 plot without y-axis tick labels
  \nextgroupplot[
    title={(e) Eigenvalues (\exref{ex:622bus})},
    xlabel={Index},
    ymode=log,
    % ymin=1e-16,
    % ymax=100,
    % yticklabels={}
  ]
  \input{5.Numerical/data/okCond_3_eig} 

    % 3 plot without y-axis tick labels
  \nextgroupplot[
    title={(f) Eigenvalues (\exref{ex: not large outlying eigs})},
    xlabel={Index},
    ymode=log,
    % ymin=1e-10,
    % ymax=1e2,
    % yticklabels={}
  ]
  \input{5.Numerical/data/okCond_1_eig}

\end{groupplot}

\path (myplot c1r1.south west|-current bounding box.south)--
      coordinate(legendpos)
      (myplot c3r1.south east|-current bounding box.south);

\matrix[
    matrix of nodes,
    anchor=south,
    draw,
    inner sep=0.2em,
    nodes={font=\small}
] at ([yshift=-9ex]legendpos)
{
    \ref{plot:fp64} & PCG(FP64) & [5pt] &
    \ref{plot:fp64mixed} & \APCG(FP64) & [5pt] &
    \ref{plot:fp16} & \APCG(FP16) & [5pt] \\
    \ref{plot:fp64tofp32} & $ u_{z,k}=\mathrm{FP32}$ & [5pt] &
    \ref{plot:fp32tofp16} & $ u_{z,k}=\mathrm{FP16}$ & [5pt] &
    \ref{plot:residualfp64tofp32} & $ u_{r,k}=\mathrm{FP32}$ & [5pt] 
    \\
};
\end{tikzpicture}
\caption{
Convergence curves of the relative residual norm for PCG(FP64), \APCG(FP64), and \APCG(FP16)  applied to  \cref{ex: not large outlying eigs,ex:bcsstk,ex:622bus} in \cref{sec: matcol}.  
The red and black crosses indicate the iterations at which $u_{z,k}$, used for $\Bz_k$ and $\Bp_k$, is set to FP32 and FP16, respectively. The blue circle marks the iterations where $u_{r,k}$ is switched to FP32.
}
\label{fig:matrix collection}
\end{figure}

In \cref{fig:matrix collection}, we present the convergence curves for  PCG(FP64), \APCG(FP64), and \APCG(FP16). 
\APCG(FP64) initializes all variables in double precision and adaptively reduces their precision according to the precision selector function \cref{eq:precision_selector}. In contrast, \APCG(FP16) uses $\Bz_k$ and $\Bp_k$ in half precision throughout all iterations, while $\Br_k$ and $\Bq_k$ are initially in double  and subsequently reduced to single precision based on the same function.

Dotted lines in subfigures (a)–(c) show that \APCG(FP16) converges more slowly than PCG(FP64) for all three matrices. In contrast, \APCG(FP64), which stepwise reduces the precision of \( \Bz_k \) and \( \Bp_k \), achieves convergence nearly identical to that of PCG(FP64).
Red and black crosses mark the iterations where \( u_{z,k} \) is reduced to single and half precision, respectively. 
The precision switch thresholds in \cref{eq: criterion} trigger these reductions primarily during the linear convergence regime. As shown in \cref{sec: z and p precision}, using \( \Bz_k \) and \( \Bp_k \) in half precision in this regime has minimal effect on the convergence rate. 
Blue circles indicate where \( u_{r,k} \) is switched to single precision.
Subfigures (d)–(f) show the eigenvalue distributions of the matrices from \cref{ex:bcsstk,ex:622bus,ex: not large outlying eigs}, respectively;  none exhibit large, outlying eigenvalues.

These results demonstrate that for matrices without many large outlying eigenvalues, if the convergence is not linear in the initial iterations, it is necessary to start \( \Bp_k \) and \( \Bz_k \) in FP64 and then stepwise reduce their precision to FP32 and FP16 during the linear convergence regime. Moreover, the precision switching thresholds in \eqref{eq: criterion} work effective.

\subsection{Performance comparison}
\label{sec: performance}
We compare and analyze the performance of PCG(FP64) and \APCG(FP16) for a large-scale GPU-accelerated high-order finite element simulation.
\begin{example}
    Consider the screened Poisson equation
\begin{equation}
- \Delta u(x,y,z) + \lambda \, u(x,y,z) = f(x,y,z) \quad \text{in } \Omega=[-0.5,0.5]^3,
\label{eq:screened poisson}
\end{equation}
with homogeneous Dirichlet boundary condition and $\lambda=10^3$. The forcing function $f(x,y,z)$ is defined in \cref{eq:f}. 
The equation is discretized in a weak form  using a high-order finite element method with a polynomial degree of $N=7$ on a  $40\times 40\times 40$  tensor-product hexahedral mesh, yielding approximately $21.7 \times 10^6$ DoFs. 
To simplify the presentation of the performance analysis, we only consider the Jacobi preconditioner. 
\label{ex: jac libp}
\end{example}

Since each physical element is mapped from a reference element via a trilinear mapping, the geometry factors used in matrix-free operations are computed on the fly from the coordinates of the element vertices \cite{cao2025towards}. This only requires  \(24\) values per element, which is negligible compared to the total number of DoFs.
The experiments were performed using \textit{libParanumal} \cite{ChalmersKarakusAustinSwirydowiczWarburton2020}  on an NVIDIA A100 SXM 40GB GPU.
For half precision arrays, we use the \texttt{half2} data type and pad arrays with an additional zero entry when their length is odd. 

\begin{figure}[h!]
\centering
\begin{tikzpicture}
\begin{axis}[
    name=myplot,
    height=4.5cm,
    width=6cm,
    xlabel={Iteration},
    ylabel={Relative residual norm},
    ymode=log,
    ymin=1e-11,
    ymax=100,
    grid=both,
]
  \input{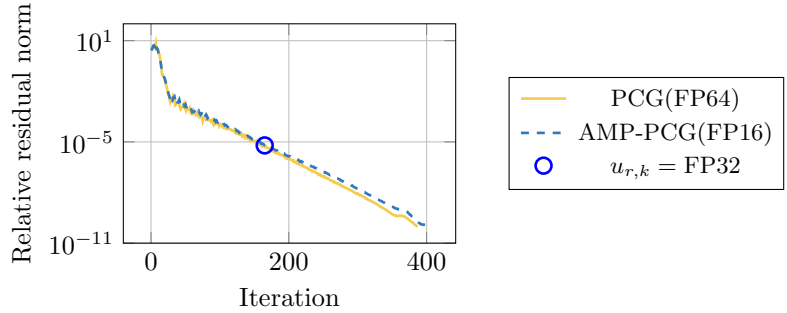}
\end{axis}

\path (myplot.south west |- current bounding box.south) --
      coordinate(legendpos)
      (myplot.south east |- current bounding box.south);

\matrix[
    matrix of nodes,
    anchor=west,
    draw,
    inner sep=0.2em,
    nodes={font=\small}
] at ([xshift=2em]myplot.east) {
    \ref{plot:fp64} &  PCG(FP64)\\
    \ref{plot:fp64mixed} & \APCG(FP16)\\
    % \ref{plot:fp64tofp32} & $ u_{z,k}=\mathrm{FP32}$ \\
    % \ref{plot:fp32tofp16} & $ u_{z,k}=\mathrm{FP16}$ \\
    \ref{plot:residualfp64tofp32} & $ u_{r,k}=\mathrm{FP32}$ \\
    \\
};
\end{tikzpicture}
\caption{%
Convergence curves of the relative residual norm for  PCG(FP64) and \APCG(FP16), applied to the screened Poisson equation with high-order finite element discretization and Jacobi preconditioner, as discussed in \cref{ex: jac libp}. $\Bz_k$ and $\Bp_k$ are in half precision throughout the iterations. The blue circles mark when $u_{r,k}$ is switched from double precision to single precision.
}
\label{fig: performance convg}
\end{figure}

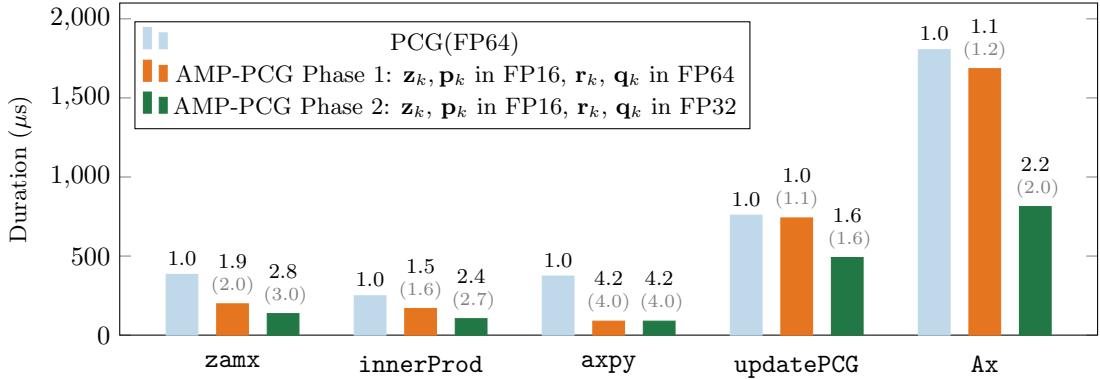
\begin{figure}[htbp]
    \centering
  \begin{tikzpicture}
\begin{axis}[
    ybar=0.25cm,
    x=2.5cm,
    bar width=12pt,
    width=15cm,
    height=6cm,
    ymin=0,
    ymax=2100,
    enlarge x limits=0.15,
    symbolic x coords={zamx, innerProd, axpy, updatePCG, Ax},
    xtick=data,
    xticklabels={\texttt{zamx}, \texttt{innerProd}, \texttt{axpy}, \texttt{updatePCG}, \texttt{Ax}},
xticklabel style={font=\ttfamily},
    ylabel={Duration (\(\mu\)s)},
    ylabel style={font=\small},
    major x tick style=transparent,
    legend style={
        at={(0.33, 0.62)},
        anchor=south,
        legend columns=1,
        font=\small,
    },
    nodes near coords,
    every node near coord/.append style={
        font=\footnotesize,
        text=black
    },
    nodes near coords style={
        /pgf/number format/fixed,
        /pgf/number format/precision=0
    }
]
   \addplot[draw=mylightblue, fill=mylightblue,
     point meta=explicit symbolic,
    nodes near coords style={
        /pgf/number format/fixed,
        /pgf/number format/precision=2,
        text=black
    }
    ] 
 coordinates { 
(zamx, 383.975) [1.0]
 (innerProd, 248.695)  [1.0]
 (axpy, 373.464)  [1.0]
 (updatePCG, 758.335)  [1.0]
 (Ax, 1804.797)  [1.0]
 }; 
   \addplot[draw=vtorange, fill=vtorange, 
       point meta=explicit symbolic,
    nodes near coords style={
        /pgf/number format/fixed,
        /pgf/number format/precision=2,
        text=black,
            yshift=10pt  % <-- increase for more spacing
    }
    ] 
 coordinates { 
(zamx, 197.925) [1.9]
 (innerProd, 168.583)  [1.5]
 (axpy, 88.380)  [4.2]
 (updatePCG, 741.092)  [1.0]
 (Ax, 1685.107)  [1.1]
 }; 
\node[font=\scriptsize, anchor=south, text=gray] at (axis cs:zamx, 197.925) {(2.0)};
\node[font=\scriptsize, anchor=south, text=gray] at (axis cs:innerProd, 168.583) {(1.6)};
\node[font=\scriptsize, anchor=south, text=gray ] at (axis cs:axpy, 88.380) {(4.0)};
\node[font=\scriptsize, anchor=south, text=gray ] at (axis cs:updatePCG, 741.092) {(1.1)};
\node[font=\scriptsize, anchor=south, text=gray ] at (axis cs:Ax, 1685.107) {(1.2)};

   \addplot[draw=mygreen, fill=mygreen,
       point meta=explicit symbolic,
    nodes near coords style={
        /pgf/number format/fixed,
        /pgf/number format/precision=2,
        text=black,
            yshift=10pt  % <-- increase for more spacing
    }] 
 coordinates { 
(zamx, 135.986) [2.8]
 (innerProd, 104.076)  [2.4]
 (axpy, 88.380)  [4.2]
 (updatePCG, 490.667)  [1.6]
 (Ax, 812.908)  [2.2]
 }; 
 \node[font=\scriptsize, anchor=south, text=gray, xshift=20pt] at (axis cs:zamx, 135.925) {(3.0)};
\node[font=\scriptsize, anchor=south, text=gray, xshift=20pt] at (axis cs:innerProd, 104.583) {(2.7)};
\node[font=\scriptsize, anchor=south, text=gray , xshift=20pt] at (axis cs:axpy, 88.380) {(4.0)};
\node[font=\scriptsize, anchor=south, text=gray , xshift=20pt] at (axis cs:updatePCG, 490.092) {(1.6)};
\node[font=\scriptsize, anchor=south, text=gray , xshift=20pt] at (axis cs:Ax, 812.107) {(2.0)};
\legend{
  {PCG(FP64)},
  {\APCG{} Phase 1: $\mathbf{z}_k, \mathbf{p}_k$ in FP16, $\mathbf{r}_k$, $\mathbf{q}_k$ in FP64},
  {\APCG{} Phase 2: $\mathbf{z}_k$, $\mathbf{p}_k$ in FP16, $\mathbf{r}_k$, $\mathbf{q}_k$ in FP32}
}
    \end{axis}
\end{tikzpicture}
\caption{
Runtime of the main kernels in PCG(FP64) and \APCG(FP16) for \cref{ex: jac libp}, measured on an NVIDIA A100 SXM 40GB GPU using \textit{Nsight Compute}.
Each bar of phase 1 and 2 is annotated with two values: the actual speedup (top line, black) relative to PCG(FP64), and the theoretical speedup (bottom line, in parentheses), estimated based on data movement detailed in \cref{tab:pcg_cost}.
The convergence behavior is shown in \cref{fig: performance convg}.
With a stopping criterion of $\|\Br_k\| \leq 10^{-10} \|\Bb\|$, \APCG(FP16) achieves a $1.63\times$ speedup over the PCG(FP64) baseline.
}
    \label{fig:performance_runtime}
\end{figure}

In \cref{fig: performance convg}, we present the convergence curves for PCG(FP64) and \APCG(FP16), where the precision  \( u_{z,k} \) is fixed at FP16 and the precision  \( u_{r,k} \) is adaptively selected based on \cref{eq:linear criterion}. 
The convergence of \APCG(FP16) closely follows that of the PCG(FP64) in this example.
\APCG(FP16) process consists of two phases:  
(1) \(\Bz_k\) and \(\Bp_k\) are stored in FP16, while all other vectors remain in FP64;  
(2) \(\Bz_k\) and \(\Bp_k\) are in FP16, \(\Br_k\) and \(\Bq_k\) are in FP32, and the remaining vectors are in FP64.
In \cref{fig:performance_runtime}, we show the runtime comparison of PCG(FP64) and the two phases of \APCG(FP16) across five main kernels: \texttt{zamx}, \texttt{innerProd}, \texttt{axpy}, \texttt{Ax}, and \texttt{updatePCG}, measured with \textit{Nsight Compute}. 
Among these, 
the \texttt{Ax} kernel is compute-bound and is expected to achieve a \(2\times\) speedup since FP32 offers twice the FLOP throughput compared to FP64 on the A100.
The other kernels are memory-bound, so their expected speedup can be estimated based on the amount of data movement involved, as summarized in \cref{tab:pcg_cost}.
Each bar in \cref{fig:performance_runtime} is annotated with two values: the actual speedup (top line) and the theoretical speedup calculated using \cref{tab:pcg_cost} (bottom line), both relative to the double precision implementation. 
The actual and theoretical speedups are closely aligned, with significant gains observed in Phase 2 of \APCG{}. With a stopping criterion of \( \| \Br_k \| \leq 10^{-10} \| \Bb \| \), \APCG(FP16) achieves an overall speedup of \( 1.63 \times \) compared to PCG(FP64).

\begin{table}[htbp]
  \centering
  \setlength{\tabcolsep}{8pt}
  \renewcommand{\arraystretch}{1.2} % increases row height by 1.5×
  \caption{Data movement (read and write) per degree of freedom in PCG(FP64) and \APCG(FP16), measured in bytes. Here, $\BM$ denotes the Jacobi preconditioner. $N_{\ell}$ is the number of local DoFs and $N_s$ is the number of global DoFs.}
  \label{tab:pcg_cost}
  \begin{tabular}{|l|l|c|c|c|}
    \hline
    \multirow{ 3}{*}{Kernel} & {Phase} & {PCG(FP64)} & {Phase 1} & {Phase 2} \\
    \cline{2-5}
    & \multirow{ 2}{*}{Precision} & \multirow{ 2}{*}{All in FP64} & {$u_{z,k}=$ FP16} & $u_{z,k}=$ FP16 \\
    & & & $u_{r,k}=$ FP64 &  $u_{r,k}=$ FP32 \\
    \hline
    \texttt{zamx} & $\Bz = \BM \Br$ & 24 & 12 & 8 \\       
    \texttt{innerProd} &$\gamma = \Bz^T\Br$ &  16 & 10 & 6 \\
     \texttt{axpy} &$\Bp = \Bz + \beta\Bp$ & 24 & 6 & 6 \\
    \texttt{updatePCG} & $\Br = \Br - \alpha\Bq$, $\Bx = \Bx + \alpha\Bp$, $\|\Br\|$ &  48 & 42 & 30 \\
    \texttt{Ax} &$\Bq = \BA\Bp$, $\mu = \Bp^T\Bq$ & $\approx 28\,{N_{\ell}}/{N_s}$ & $\approx 24\,{N_{\ell}}/{N_s}$ & $\approx 14\,{N_{\ell}}/{N_s}$ \\
     \hline
  \end{tabular}
\end{table}

\subsection{Limitations of the precision switch criterion}
\label{sec: limitation}

 In the previous subsections, we  verified the effectiveness of \APCG{} from the perspectives of attainable accuracy, convergence rate, and GPU performance across various examples, where the assumptions for the convergence are satisfied.
In \cref{sec: criterion rq}, we assume the residual decays rapidly after $d$ iterations so that \cref{eq: res norm decay} holds. This condition ensures that $\eta_k$ provides a good approximation to $ \widehat{\eta}_m$ for $m \gg k$.  
\cref{thm: 3.1} guarantees that the linear convergence rate of the perturbed $\Bz_k$ is close to that of the exact arithmetic PCG. However, there is no theoretical guarantee on the degradation of the convergence rate when the residual $\Br_k$ exhibits oscillatory behavior. In this subsection, we present several examples where these assumptions are not satisfied to illustrate limitations of the precision selector  \cref{eq:precision_selector}.

\subsubsection{Choice of $d$ in the attainable accuracy indicator}
\label{sec: choice of d}
The effectiveness of the attainable accuracy indicator $\eta_k$ depends on the delay parameter $d$, which should be large enough such that the inequality \cref{eq: res norm decay} holds. 
In previous examples, we set \( d = 10 \). However, if the convergence is slow or oscillatory, a larger value of \( d \) may be required to obtain a reliable estimate of the final attainable accuracy.
To illustrate this, we consider the matrix \texttt{nos6.mat}  from the SuiteSparse Matrix Collection \cite{davis2011university}. 
\begin{example}
The matrix \texttt{nos6.mat}, of size \(675 \times 675\), is scaled on both sides by the square root of the inverse of its diagonal entries. The condition number after scaling is approximately \(3.49 \times 10^{6}\).
    \label{ex: limit res}
\end{example}

We set the stopping tolerance $\tol=10^{-6}$. 
In \cref{fig:limit res}, we show the convergence curves for both the updated residual (left) and the true residual (right) of PCG and \APCG($\urk$) .
In these experiments, we do not terminate the iterations based on the stopping condition in line 18 in \cref{alg:MPDS-PCG}, allowing us to observe the final attainable accuracy achieved by the algorithm. 
Markers indicate the iteration at which the precision of $\Br_k$ and $\Bq_k$ is switched from FP64 to FP32 base on the precision selector \cref{eq:precision_selector}.
All other vectors and operations are kept in FP64.
In this case, the residual norm initially decreases and then increases, so the condition \cref{eq: res norm decay}  does not hold for \( d = 10 \). As a result, \( \eta_k \) fails to provide an upper bound for the attainable accuracy, and the final relative true residual norm exceeds $\tol = 10^{-6}$.
To improve the accuracy of \( \eta_k \) as an estimator, a larger delay parameter, at least $d=50$, would be needed.
The difficulty of selecting a proper $d$ also arises in error estimation for CG  using the delay strategy \cite{hestenes1952methods, arioli2004stopping}, and it is also an open research area.

\begin{figure}[htbp]
\centering
\begin{tikzpicture}
\begin{groupplot}[
    group style={
        group name=myplot,
        group size=2 by 1,
        horizontal sep=20pt,
        vertical sep=2cm
    },
    height=4.5cm,
    width=5.5cm
]
  % First plot with y-axis labels
  \nextgroupplot[
  title={(a) Updated residual},
    xlabel={Iteration},
    ylabel={Relative residual norm},
    ymode=log,
    ymin=1e-15,
    ymax=20,
    grid=both,
        ytick={1e0,  1e-6, 1e-12},
  ]
  \input{5.Numerical/data/resPrec_limit_1}
  
  % Second plot without y-axis tick labels
  \nextgroupplot[
    title={(b) True residual},
    xlabel={Iteration},
    ymode=log,
    ymin=1e-15,
    ymax=20,
    grid=both,
        ytick={1e0,  1e-6, 1e-12},
    yticklabels={},
  ]
  \input{5.Numerical/data/resPrec_limit_2}
 
\end{groupplot}

\path (myplot c1r1.south west|-current bounding box.south)--
      coordinate(legendpos)
      (myplot c2r1.south east|-current bounding box.south);

\matrix[
    matrix of nodes,
    anchor=west, % align the left edge of the legend with the anchor point
    draw,
    inner sep=0.2em,
    nodes={font=\small}
] at ([xshift=2em]myplot c2r1.east) % shift a bit to the right
{
    \ref{plot:fp64} & PCG(FP64) \\
    \ref{plot:resPrec1e6_limit} & \APCG($\urk$)  \\
    \ref{plot:resPrec1e6_ind_limit} & $\Weta_k$   \\
    \ref{plot:residual2_limit} & Switch point  \\
};
\end{tikzpicture}
\caption{Convergence curves of (a) relative updated residual norm and (b) relative true residual norm for  PCG(FP64), \APCG($\urk$)  with a stopping tolerance  $10^{-6}$  applied to \cref{ex: limit res}. In \APCG($\urk$) , the indicator \cref{eq: eta_k} with $d=10$ is used in the precision selector \cref{eq:precision_selector} to switch the precision of $\Br_k$ and $\Bq_k$ from FP64 to FP32, while all other vectors are in FP64.
$\Weta_k$ is the final attainable accuracy indicator defined in \cref{eq: Weta_k}. 
Since the residual norm decreases and then increases, the assumption \cref{eq: res norm decay} does not hold for \( d = 10 \). As a result, the final accuracy does not meet the target tolerance of \( 10^{-6} \).
 }
\label{fig:limit res}
\end{figure}

\subsubsection{Choice of precision switch thresholds $\tau_{z,s}$ and $\tau_{z,h}$}
\label{sec: choice of tauz}
In this subsection, we consider the impact of thresholds $\tau_{z,s}$ and $\tau_{z,h}$ on the convergence rate. 
Several examples in previous subsections show that applying the \APCG{} algorithm, i.e., stepwise reduction of the precision of $\Bp_k$ and $\Bz_k$ in the linear convergence regime, rarely impacts the convergence rate. 
However, when the matrix has many large outlying eigenvalues, it may require a significant number of iterations to reach the linear convergence regime. 
As a result, selecting appropriate thresholds $\tau_{z,s}$ and $\tau_{z,h}$
can be challenging.
We first consider a  $10 \times 10$ matrix to demonstrate that using the thresholds  \cref{eq: criterion} to switch the precision of $\Bz_k$ and $\Bp_k$ in \APCG{} leads to slow convergence relative to PCG(FP64)
We then present another example illustrating how the choice of $\tau_{z,s}$ impacts the convergence of \APCG{}.
In  \cite{greenbaum1992predicting, carson2024towards}, it is shown that the convergence rate of PCG in finite precision can degrade when the system matrix has multiple large outlying eigenvalues. We consider a similar matrix in the following example.
\begin{example}
Consider an SPD matrix with eigenvalues $\lambda_i\, (i = 1, \dots, n)$. Its largest eigenvalue is $\lambda_n = 1$,  while many eigenvalues are clustered near $\lambda_1 = 10^{-3}$, defined as 
\begin{equation}
    \lambda_i = 10^{-3} + \frac{i-1}{n-1} \, \rho^{\,n-i},\, i=1,\dots,n-1,
\label{eq:limitation}
\end{equation}
where $n = 10$ and $\rho = 1/4$. $\Bb$ is set to a vector of all ones.
\label{ex: limitation 1}
\end{example}

\begin{figure}[htbp]
\centering
\begin{tikzpicture}
\begin{groupplot}[
    group style={
        group name=myplot,
        group size=1 by 1, % Only one plot
        horizontal sep=20pt,
        vertical sep=2cm
    },
    height=4.5cm,
    width=5.5cm,
]
  \nextgroupplot[
    title={Residual},
    xlabel={Iteration},
    ylabel={Relative residual norm},
    ymode=log,
    ymin=1e-15,
    ymax=20,
    grid=both,
  ]
  \input{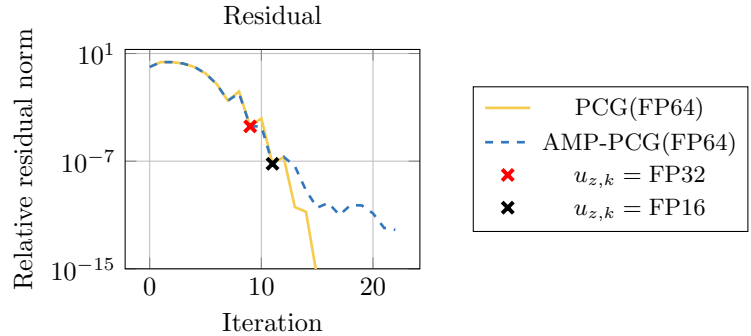}
\end{groupplot}

% Fix legend positioning — only one plot, so use c1r1
\path (myplot c1r1.south west |- current bounding box.south) --
      coordinate(legendpos)
      (myplot c1r1.south east |- current bounding box.south);

% Clean 2-column legend matrix
\matrix[
    matrix of nodes,
    anchor=west,
    draw,
    inner sep=0.2em,
    nodes={font=\small}
] at ([xshift=2em]myplot c1r1.east) {
    \ref{plot:fp64} & PCG(FP64) \\
    \ref{plot:fp64mixed} & \APCG(FP64) \\
    \ref{plot:fp64tofp32} & $u_{z,k} = \mathrm{FP32}$ \\
    \ref{plot:fp32tofp16} & $u_{z,k} = \mathrm{FP16}$ \\
};

\end{tikzpicture}
\caption{Convergence curves of the relative residual norm for PCG(FP64) and \APCG(FP64), applied to \cref{ex: limitation 1}. The red and black crosses indicate the iterations at which $u_{z,k}$, used for $\Bz_k$ and $\Bp_k$, is set to FP32 and FP16 in \APCG(FP64), respectively.
}
\label{fig:limit small convg}
\end{figure}

We compare PCG(FP64)  with \APCG(FP64), where $\Bz_k$ and $\Bp_k$ are initially in double precision and stepwise reduced to single and half precision based on the thresholds in \cref{eq: criterion}. 
To isolate the impact of inexactness in $\Bz_k$ and $\Bp_k$, all other vectors are maintained in double precision.
Figure~\ref{fig:limit small convg} shows the convergence of both variants. 
In \APCG(FP64), the precision of $\Bz_k$ and $\Bp_k$ is reduced to single precision at iteration 9 (red cross), and further reduced to half precision at iteration 11 (black cross).
This progressive reduction in precision results in a few additional iterations compared to the PCG(FP64).

\begin{figure}[htbp]
    \centering
    \begin{subfigure}{0.49\textwidth}
        \centering
        \includegraphics[width=\linewidth]{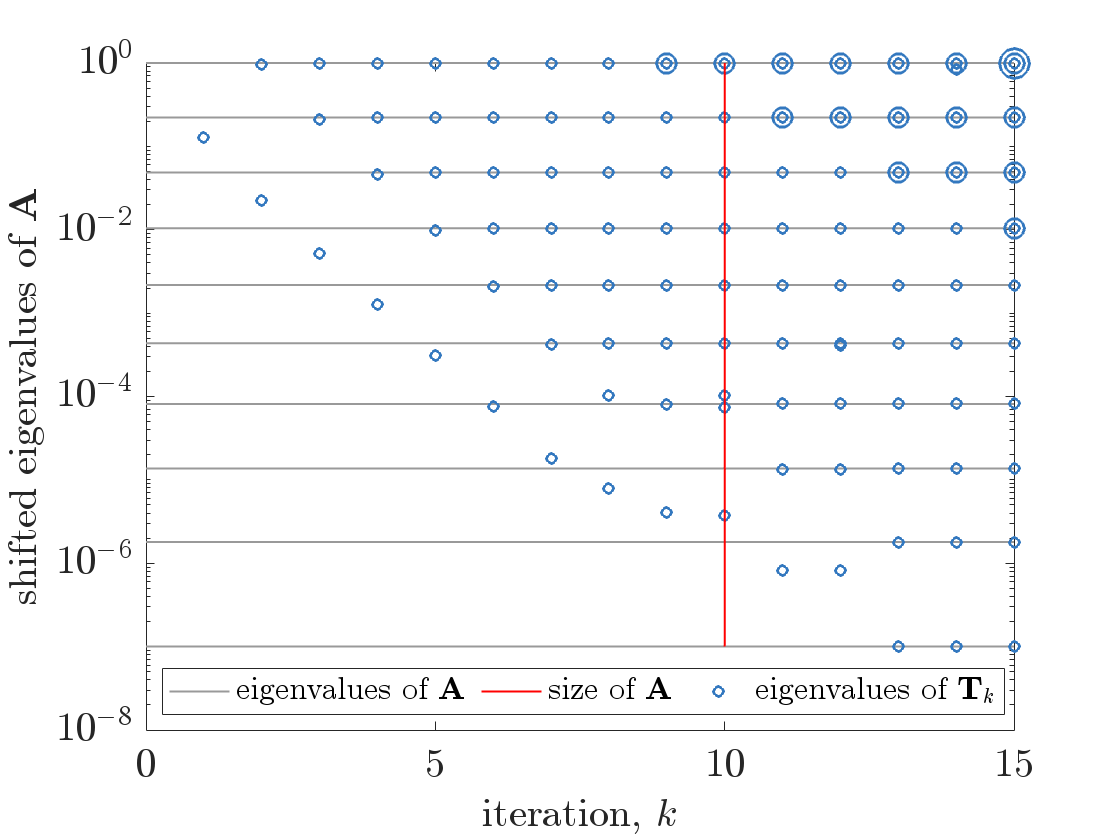}
        \caption{Ritz values from PCG(FP64).}
        \label{fig:top}
    \end{subfigure}
    \hfill
    \begin{subfigure}{0.49\textwidth}
        \centering
        \includegraphics[width=\linewidth]{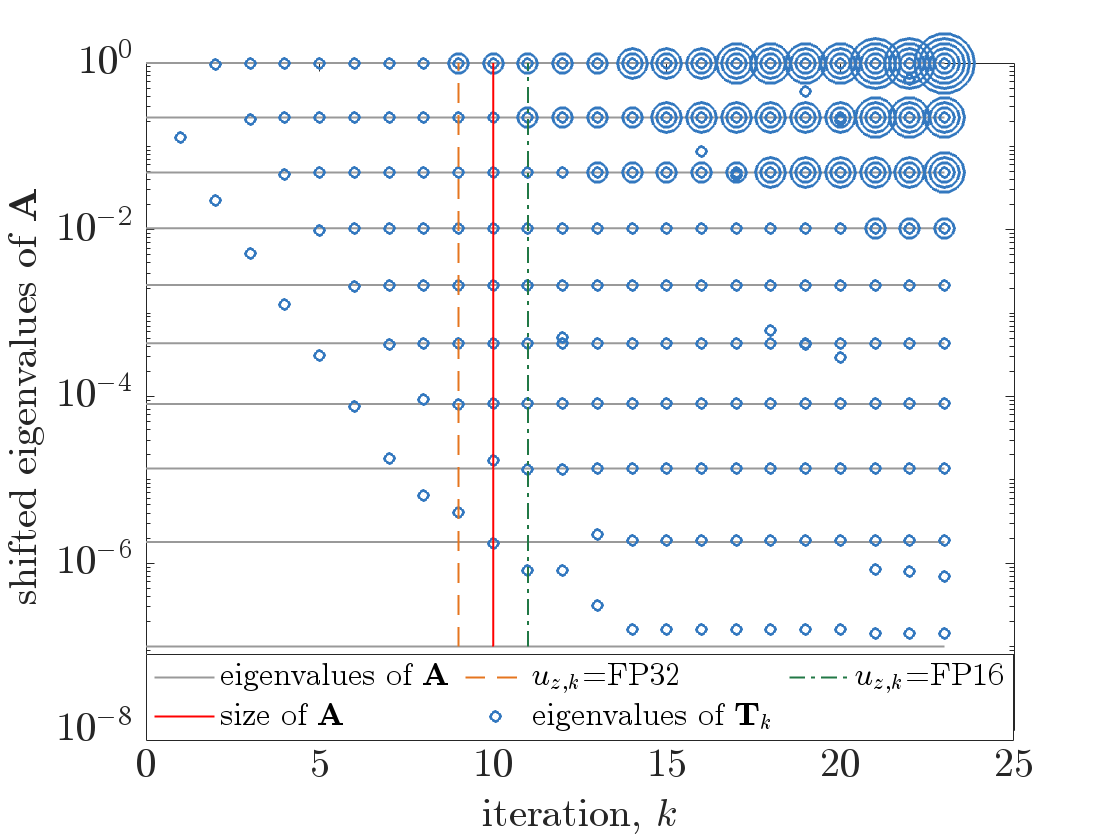}
        \caption{Ritz values from \APCG(FP64).}
        \label{fig:bottom}
    \end{subfigure}
    \caption{
Comparison of Ritz values generated by PCG(FP64) and \APCG(FP64) for \cref{ex: limitation 1}. The corresponding convergence curves are shown in \cref{fig:limit small convg}. For better visualization, the $y$-axis shows 
eigenvalues shifted by \( - 10^{-3} + 10^{-7} \), which maps the smallest eigenvalue to $10^{-7}$.
% shifted eigenvalues, defined as \( \lambda_i - 10^{-3} + 10^{-7} \).
In \APCG(FP64), the precision of \( \Bz_k \) and \( \Bp_k \) is reduced from FP64 to FP32 at the 9\textsuperscript{th} iteration (dashed orange line), and further to FP16 at the 11\textsuperscript{th} iteration (dash-dotted green line), while all other vectors remain in FP64.
Subfigure (a) shows that PCG(FP64) approximates the largest eigenvalue three times in the final iteration, whereas subfigure (b) shows that \APCG(FP64) approximates the largest eigenvalues six times. This repeated approximation of large, outlying eigenvalues contributes to the slower convergence observed in \APCG(FP64).
    }
    \label{fig:eigen_convg}
\end{figure}

In exact arithmetic, CG should converge in at most 10 iterations since $\BA\in\mathbb{R}^{10\times 10}$. However, both PCG(FP64) and \APCG(FP64) require more than 10 iterations to converge. This behavior is explained in \cite{greenbaum1992predicting, carson2024towards}, where results indicate that in finite precision arithmetic, multiple approximations to the larger eigenvalues may appear before any accurate approximation to the smaller, clustered eigenvalues is obtained.
\cref{fig:eigen_convg} illustrates the Ritz values, i.e. the  eigenvalues of the Lanczos tridiagonal matrices $\BT_k$, as the iteration progresses. 
The $x$-axis represents the iteration, while the $y$-axis shows the eigenvalues shifted by  $- 10^{-3} + 10^{-7}$, which maps the smallest eigenvalue to $10^{-7}$, allowing for clearer visualization.
 The gray horizontal lines represent eigenvalues of $\BA$. 
Circles denote the Ritz values at each iteration and circles at the same position represent repeated approximation to the same eigenvalue. 
In the PCG(FP64), the largest eigenvalue is approximated three times. In contrast, \APCG(FP64) detects large eigenvalues even more frequently; the largest eigenvalue is approximated six times, and the next two largest eigenvalues have multiplicity of four in the final iteration. This repeated discovery of large, outlying eigenvalues contributes to slower convergence, as it leads to the loss of linear independence of the Krylov basis vectors.

Next, we consider a test case solved by \APCG{} using two different thresholds,  $\tau_{z,s} = 10^{-4}$ (as in \cref{eq: criterion}) and  $\tau_{z,s} = 10^{-8}$. We set $\tau_{z,h} = 0$ (meaning conversion to half precision is never triggered), and all other vectors are represented in double precision.
\begin{example}
Consider a matrix $\BA$ with the eigenvalues defined by \cref{eq:limitation} with $n = 50$ and $\rho = 0.9$.
\label{ex: limitation 2}
\end{example}

\begin{figure}[htbp]
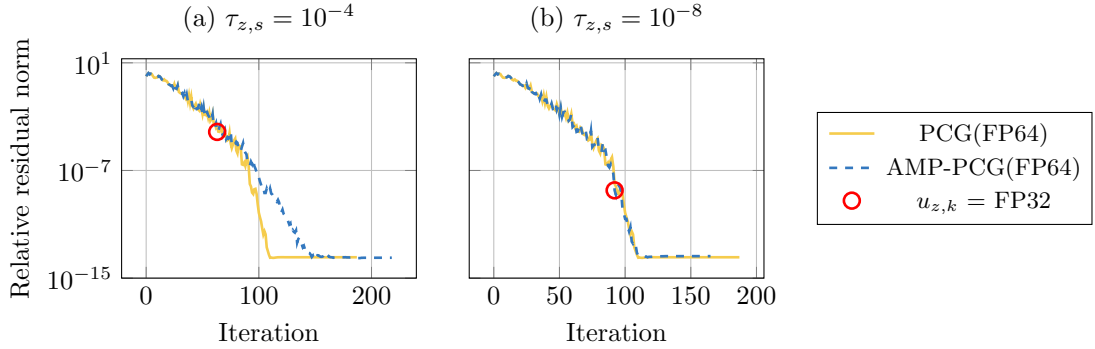

\begin{tikzpicture}
\begin{groupplot}[
    group style={
        group name=myplot,
        group size=3 by 1,
        horizontal sep=20pt,
        vertical sep=2cm
    },
    height=4.5cm,
    width=5.5cm
]
  % First plot with y-axis labels
  \nextgroupplot[
  title={(a) $\tau_{z,s} = 10^{-4}$},
    xlabel={Iteration},
    ylabel={Relative residual norm},
    ymode=log,
    ymin=1e-15,
    ymax=20,
    grid=both,
  ]
  \input{5.Numerical/data/limitation1}
  
  % Second plot without y-axis tick labels
  \nextgroupplot[
    title={(b) $\tau_{z,s} = 10^{-8}$},
    xlabel={Iteration},
    ymode=log,
    ymin=1e-15,
    ymax=20,
    yticklabels={},
    grid=both,
  ]
  \input{5.Numerical/data/limitation2}
  
\end{groupplot}

\path (myplot c1r1.south west|-current bounding box.south)--
      coordinate(legendpos)
      (myplot c2r1.south east|-current bounding box.south);

\matrix[
    matrix of nodes,
    anchor=west, % align the left edge of the legend with the anchor point
    draw,
    inner sep=0.2em,
    nodes={font=\small}
] at ([xshift=2em]myplot c2r1.east) % shift a bit to the right
{
    \ref{plot:fp64} & PCG(FP64) \\
    \ref{plot:fp64mixed} & \APCG(FP64) \\
    % \ref{plot:limfp16} &  DSCG with $\Bz_k,\,\Bp_k$ in FP32 & [5pt] \\
    % \ref{plot:limfp32} & DSCG with $\Bz_k,\,\Bp_k$ in FP16& [5pt] &
    \ref{plot:limfp64tofp32} & $u_{z,k}$ = FP32 \\
    % \ref{plot:limfp32tofp16} & set $u_{z,k}$ = FP16 & [5pt] \\
};
\end{tikzpicture}
\caption{
Convergence curves of the relative true residual norm of  PCG(FP64) and \APCG(FP64) for \cref{ex: limitation 2}. 
In  \APCG(FP64), the precision of $\Bz_k$ and $\Bp_k$ starts in FP64 and then are switched to FP32 at red circles. 
$\tau_{z,s}$ is used in \cref{eq:precision_selector} to determine when to switch precision of $\Bz_k$ and $\Bp_k$.  
All other vectors are in FP64.
}
\label{fig:carson example}
\end{figure}

In \cref{fig:carson example}, we present the convergence curves of the relative residual norm for  PCG in FP64 and  \APCG(FP64), where $\Bz_k$ and $\Bp_k$ start in FP64 and are switched to FP32 based on $\tau_{z,s}$, with red circles indicating the switch points.
In subfigure (a), with $\tau_{z,s} = 10^{-4}$, the threshold used in \cref{eq: criterion} and previous examples, the precision switch occurs during an oscillatory phase of the convergence. 
After switching to FP32, \APCG(FP64) exhibits nearly linear convergence, whereas PCG(FP64) continues to show superlinear convergence. 
In contrast, in subfigure (b), with $\tau_{z,s} = 10^{-8}$, the switch occurs closer to the linear convergence regime, and the convergence rates of PCG(FP64) and \APCG(FP64) are almost identical.

This example underscores the sensitivity of the convergence rate of \APCG{} to the precision switch threshold. Specifically, when a matrix has multiple large, outlying eigenvalues, those eigenvalues tend to be approximated frequently, requiring a substantial number of iterations before the linear convergence regime is reached. Consequently, selecting the appropriate iteration to switch precision is nontrivial. The thresholds in \cref{eq: criterion} are empirical and do not perform well in this case.
Furthermore, if \( n = 100 \) is used in \cref{ex: limitation 2}, the convergence does not exhibit a clearly identifiable linear regime. In this case, reducing the precision of \( \Bz_k \) and \( \Bp_k \) at any iteration leads to slower convergence compared to PCG(FP64). Hence, for AMP-PCG, preconditioning should aim to avoid or reduce
large, outlying eigenvalues of the
preconditioned system.
The impact of roundoff errors in PCG has been widely studied in the context of finite precision  \cite{greenbaum1997iterative, meurant2006lanczos, carson2024towards}, and these insights are directly relevant for understanding the convergence behavior of \APCG{}.

\section{Summary}

We have presented an adaptive mixed precision and dynamically scaled preconditioned conjugate gradient algorithm for solving sparse linear systems. As demonstrated in the motivating example in \cref{sec:motivation}, naively implementing PCG in low precision faces three major challenges: reduced attainable accuracy, delayed convergence, and numerical underflow due to the limited dynamic range of half precision.
The proposed \APCG{} algorithm addresses these issues by adaptively selecting the precision of  vectors and dynamically scaling the right-hand side vector in the preconditioning step to ensure all vectors remain representable in low precision. An attainable accuracy indicator is introduced to guide the precision switching of $\Br_k$ and $\Bq_k$ from double to single precision. Additionally, we propose an empirical criterion for adjusting the precision of  $\Bz_k$ and $\Bp_k$: if convergence is nearly linear, both vectors remain in half precision throughout the iterations; if convergence is superlinear, they begin in double precision and switch to lower precision once the linear convergence is reached.

Our experiments on a broad set of sparse linear systems demonstrate the strong potential of the proposed algorithm. 
We show that across a wide variety of problems, the convergence rate remains unaffected by adaptive mixed precision strategies, with final solution accuracy comparable to uniform double precision implementations. 
We also analyze the performance of a large-scale GPU-accelerated high-order finite element discretization with a Jacobi preconditioner, in which $\Bz_k$ and $\Bp_k$ are maintained in half precision throughout the iterations. The vectors $\Br_k$ and $\Bq_k$ initially use double precision and are later reduced to single precision. This use of low precision significantly reduces data movement and floating-point operations in key kernels, achieving an overall speedup of about $1.63$ on an A100 GPU while maintaining the accuracy of the computed solution.

The potential shortcomings of the precision selector function \cref{eq:precision_selector} are discussed in \cref{sec: limitation}. For matrices with many large, outlying eigenvalues, selecting an appropriate threshold for switching the precision of $\Bz_k$ and $\Bp_k$ can be challenging. An improper switch may degrade the convergence rate. This highlights the need for either a more robust precision-switching criterion or the use of an effective preconditioner for such matrices.
Another important consideration is the trade-off between iteration count and per-iteration cost. Although \APCG{} may require more iterations to converge, each iteration can be significantly cheaper. Balancing convergence rate and computational cost per iteration remains an interesting direction for further investigation.
While our performance study focuses on linear systems arising from high-order finite element discretizations with a simple Jacobi preconditioner, greater speedups are expected when using more effective preconditioners, such as \( p \)-multigrid. In future work, we plan to optimize low precision implementations of multigrid preconditioners and apply \APCG{} to exascale, real-world, time-dependent problems.

\appendix

\section{Proof of \cref{thm: dspcg}}
\label{apx: dspcg}
\begin{proof}
All variables produced in the \DPCG{} process are denoted by hats. We proceed it by induction.
Starting with the first iteration ($k=1$) 
\[
\Hz_0 = \omega_0 \BM^{-1}\Hr_0 = \omega_0 \BM^{-1}\Br_0 = \omega_0 \Bz_0, \quad \Hrho_0 = \omega_0 \rho_0, \quad \Hp_0 = \Hz_0 =  \omega_0 \Bp_0, \quad  \Hq_0 = \omega_0 \Bq_0.
\]
Therefore, we obtain $\Hgamma_0 = \omega_0^2 \gamma_0, \,\Halpha_0 = {\alpha_0}/{\omega_0}$, and thus 
\[
\Halpha_0 \Hp_0 = \alpha_0 \Bp_0, \quad \Halpha_0 \Hq_0 = \alpha_0 \Bq_0.
\]
Consequently, 
\[
\Hx_1 = \Bx_1 \text{ and } \, 
\Hr_1 = \Br_1.
\]
Thus, the claim holds for $k=1$.

Next, we assume that $\Hrho_{k-1} = \omega_{k-1} \rho_{k-1}$, $\Hp_{k-1} = \omega_{k-1} \Bp_{k-1}$,  $\Hx_k = \Bx_k$,  and $\Hr_k = \Br_k$ for some positive integer $k$. We aim to show that the identities also hold at the next iteration: $\Hrho_{k} = \omega_{k}\rho_{k}$, $\Hp_{k} = \omega_{k}\Bp_{k}$,  $\Hx_{k+1} = \Bx_{k+1}$,  and $\Hr_{k+1} = \Br_{k+1}$.
Similar  to the initial step,
\[
\Hz_k = \omega_k\,\BM^{-1}\Hr_k = \omega_k\,\BM^{-1}\Br_k = \omega_k\,\Bz_k, \quad \Hrho_k = \omega_k\,\rho_k.
\]
Moreover, for search direction, 
\[
\Hp_k = \Hz_k + \Hbeta_{k-1}\Hp_{k-1} 
= \omega_k\Bz_k + \frac{\omega_k\rho_k}{\omega_{k-1}\rho_{k-1}}\omega_{k-1}\Bp_{k-1} 
= \omega_k \lp \Bz_k + \frac{\rho_k}{\rho_{k-1}}\Bp_{k-1} \rp = \omega_k \Bp_{k}.
\]
By multiplying $\BA$ to both sides, we obtain $\Hq_k = \omega_k\,\Bq_k$.
Then,
\[
\Hgamma_k = \Hq_k^T\Hp_k = (\omega_k\,\Bq_k)^T (\omega_k\,\Bp_k) = \omega_k^2\,\gamma_k.
\]
Thus, the step size is
\[
\Halpha_k = \frac{\Hrho_k}{\Hgamma_k}  = \frac{\rho_k}{\omega_k\,\gamma_k} = \frac{\alpha_k}{\omega_k},
\]
with $\alpha_k = \rho_k/\gamma_k$. Hence, the solution update is
\[
\Hx_{k+1} = \Hx_k + \Halpha_k\,\Hp_k = \Bx_k + \frac{\alpha_k}{\omega_k}(\omega_k\,\Bp_k) = \Bx_k + \alpha_k\,\Bp_k = \Bx_{k+1},
\]
and the residual update becomes
\[
\Hr_{k+1} = \Hr_k - \Halpha_k\,\Hq_k = \Br_k - \frac{\alpha_k}{\omega_k}(\omega_k\,\Bq_k) = \Br_k - \alpha_k\,\Bq_k = \Br_{k+1}.
\]
By induction, the result holds for all iterations $k$. That is, for every $k$, 
$\Hx_{k} = \Bx_{k}$  and $\Hr_{k} = \Br_{k}$.
\end{proof}

% \section*{Acknowledgments}
% We would like to acknowledge the assistance of volunteers in putting
% together this example manuscript and supplement.

\bibliographystyle{siamplain}
\bibliography{ref}

\end{document}